\documentclass[final,leqno]{siamltex}

\usepackage{hyperref}
\usepackage{amsmath,amsfonts}
\usepackage{xcolor,graphicx}
\usepackage{multirow}
\usepackage{enumitem}
\setlist[enumerate]{label*=\arabic*.}

\newcommand{\abs}[1]{\left\vert #1\right\vert}
\newcommand{\adim}[1]{\tilde{#1}}

\newcommand{\ext}[1]{{#1}'}
\newcommand{\norm}[2]{\left\Vert #1\right\Vert_{#2}}
\newcommand{\R}{\mathbb{R}}
\newcommand{\unit}[1]{~\mathrm{#1}}
\newcommand{\rhomax}{\rho_{\textup{max}}}
\newcommand{\tevac}{t_{\textup{evac}}}

\DeclareMathOperator*{\argmax}{arg\,max}
\DeclareMathOperator{\Div}{div}

\graphicspath{{./figure/}}

\newtheorem{remark}[theorem]{Remark}

\title{Modeling rationality to control self-organization of crowds: An environmental approach\thanks{This research has received funding from the European Union FP7 under grant No.\ 257462 HYCON2 Network of Excellence. It was also partially funded by the Google Research Award ``Multipopulation Models for Vehicular Traffic and Pedestrians'', 2012-2013.}}

\author{Emiliano Cristiani\thanks{Istituto per le Applicazioni del Calcolo ``M. Picone'', Consiglio Nazionale delle Ricerche, Rome, Italy.
\texttt{\{e.cristiani, f.priuli, a.tosin\}@iac.cnr.it}} 
\and	
Fabio S. Priuli\footnotemark[2] 
\and	
Andrea Tosin\footnotemark[2]	
}
		
\begin{document}
\maketitle

\begin{abstract}
In this paper we propose a classification of crowd models in built environments based on the assumed pedestrian ability to foresee the movements of other walkers. At the same time, we introduce a new family of macroscopic models, which make it possible to tune the degree of predictiveness (i.e., rationality) of the individuals. By means of these models we describe both the natural behavior of pedestrians, i.e., their expected behavior according to their real limited predictive ability, and a target behavior, i.e., a particularly efficient behavior one would like them to assume (for, e.g., logistic or safety reasons). Then we tackle a challenging shape optimization problem, which consists in controlling the environment in such a way that the natural behavior is as close as possible to the target one, thereby inducing pedestrians to behave more rationally than what they would naturally do. We present numerical tests which elucidate the role of rational/predictive abilities and show some promising results about the shape optimization problem.
\end{abstract}



\pagestyle{myheadings}
\thispagestyle{plain}
\markboth{E. CRISTIANI, F. S. PRIULI, AND A. TOSIN}{MODELING RATIONALITY TO CONTROL CROWDS}

\section{Introduction}
In this paper we are concerned with the modeling of a variety of human behaviors in crowds, characterized by a different degree of pedestrian ability to foresee and anticipate the movements of other walkers. In addition, we investigate the possibility to induce pedestrians to unconsciously behave more rationally than what they would consciously do.

Although the research in pedestrian modeling is relatively young, the literature in this field is already large. This is probably due to the fact that many models are inspired by the vehicular traffic literature, which has been deeply investigated for decades. Basically all kinds of models have been proposed so far: Nanoscopic, microscopic, mesoscopic, macroscopic, multiscale, cellular automata, discrete choice, in many variants: First-order, second-order, differential, non-differential, local, nonlocal, with or without contact-avoidance features. A number of reviews are already available, see, e.g.,~\cite{bellomo2011SR,duives2013TRC}. Concerning books, we mention~\cite{kachroo2008book}, which develops macroscopic models for the control of pedestrian behavior in connection with the implementation of evacuation strategies, and~\cite{cristiani2014book}, which investigates extensively the multiscale (microscopic+macroscopic) model originally proposed in~\cite{cristiani2011MMS} and presents a detailed review of the literature, including an accurate account of seminal papers.

Pedestrians are \emph{active} particles, namely they are not passively prone to external force fields like inert matter. They observe the surroundings and elaborate a walking strategy while moving. In many cases, they have a specific target and want to reach it either in minimal time or while optimizing some other personal performance criterion. The question is to what extent pedestrians are able to optimize their path and anticipate the movements of the others. Let us describe five possible behaviors in order of increasing degree of predictive ability, say of \emph{rationality}:
\paragraph{Irrational} This behavior is typical of panic situations. An irrational pedestrian is unable to plan a path, moves faster than usual, changes direction more often, and is attracted toward people who have a clear direction. Moreover, s/he is selfish, his/her interactions become physical and coordinated movements are lost.
\paragraph{Basic} A basic pedestrian decides in advance his/her path on the basis of his/her destination and of the knowledge of the environment. For instance, s/he may follow a minimal time path if s/he knows well the environment from repeated past experiences and wants to leave it as soon as possible. His/Her path might be more irregular if his/her knowledge of the environment is limited or if s/he is unhurried. Once the path has been decided, s/he can still stray from it because of the presence of other walkers, but usually diversions are only temporary. This is the effect of short-range interactions among pedestrians, who typically are repulsed by crowded regions and want to avoid mutual collisions.
\paragraph{Rational} Besides basic collision-avoidance maneuvers, a rational pedestrian \textit{continuously} reconsiders the performance of his/her path on the basis of the \emph{current} positions of all other walkers (and obstacles, if any) in the space.
\paragraph{Highly rational} 
A highly rational pedestrian is able to forecast with no errors the movements of the others at any later time, and decides his/her path on the basis of his/her destination and the position of the others at both \textit{current} and \textit{later} times. The assumption that all pedestrians in a crowd are highly rational leads to a differential game-type setup, where the collective evolution has to be computed as a whole. In this case, it is convenient to resort to well-known concepts of differential games and mean field games, such as Nash equilibrium strategies. If a Nash equilibrium is reached, no pedestrian can find a better path if the others keep their choice. Then the situation is somehow ``stuck'', since no one finds it convenient to change unilaterally his/her strategy. However, a Nash equilibrium does not represent, in general, the best choice for the players among all the possible strategies.
\paragraph{Optimal} This behavior applies only to crowds and not to single individuals. An optimal crowd moves in such a way that its global payoff (to be defined) is maximized. Individuals in the crowd can agree to penalize their own gain in favor of the group.

Note that all the behaviors described above can be actually observed. For example, pedestrians moving in a large unfamiliar environment are likely to behave basically, while pedestrians moving in a small and well-known environment are likely to behave rationally or highly rationally. Soldiers can instead behave optimally, if an optimal strategy is imposed to them. 

Most of the models available in the literature are for basic pedestrians. Let us mention, among others, the models proposed in \cite{agnelli2015M3AS,cristiani2011MMS,helbing1995PRE,xia2009PRE}. Models for panicky pedestrians are less studied; see, e.g.,~\cite{helbing2000N}. The most popular model for rational pedestrians is the one by Hughes~\cite{hughes2002TRB}, which was widely investigated in the recent years, see references in~\cite{cristiani2014book}. Regarding highly rational pedestrians, we mention the models~\cite{hoogendoorn2003OCAM,hoogendoorn2004TRBb} and, in the framework of stochastic mean-field games,~\cite{BuDFMaWo,lachapelle2011TRB}. To our knowledge, there are no models for optimal pedestrians.

A natural issue can be posed about the classification of pedestrian behaviors set forth above: Is it really worth modeling more-than-rational pedestrians? Lachapelle and Wolfram, while proposing a model for highly rational pedestrians, doubt it. Quoting from~\cite{lachapelle2011TRB}:
``\textit{We assume that [...] every individual or agent knows and anticipates the state of everybody through the distribution of the crowd. These assumptions are probably wrong in the real world (except in routine crowd situations, say the daily way towards the subway exit)}''. Indeed, the more rationality is assumed the higher the risk of simulating an unrealistic behavior. In this paper we want to adopt a different point of view: Rational models are not used to \emph{describe} but rather to \emph{control} reality. More specifically, the goal of the paper is twofold:
\paragraph{Modeling} We introduce a new family of first-order macroscopic nonlocal fun\-da\-men\-tal-diagram-free models, which allow one to turn on or off pedestrian rationality as well as to tune its degree. In this respect, it is important to note that current models for rational and highly rational pedestrians do not allow one to weaken or turn off the rationality. Conversely, in our model the rationality is modeled by an \emph{additive component} of the velocity field superimposed to an always-defined interaction component, which allows one to turn rationality on or off at will.
\paragraph{Optimization} By means of the model just sketched, for any given situation we can simulate both the \emph{natural behavior} of pedestrians, i.e., the behavior they are expected to assume in reality, and a \emph{target behavior}, i.e., a desired behavior we would like them to assume (related, e.g., to logistic or safety reasons). With both behaviors defined, we can induce the first one to align with the second one. More precisely, we can modify the environment by means of additional (controlled) obstacles in such a way that the \emph{natural behavior in this new environment} is as close as possible to the target one. As a byproduct, we recover the well-known Braess' paradox~\cite{braess2005TS,hughes2003ARFM}, which basically states that adding obstacles or constraints can improve global dynamics.

\smallskip

In order to contextualize and compare our results with the available literature, it is useful to recall some papers where a shape optimization of the environment with respect to some cost criterion has already been attempted. In the framework of microscopic models, we mention the papers~\cite{johansson2007,shukla2009}, where the optimal shape is achieved by means of a simple genetic algorithm. Numerical tests recover to some extent the Braess' paradox, producing obstacles of various shapes which improve the outflow rate of pedestrians. In the framework of macroscopic models, we mention instead paper~\cite{twarogowska2014AMM}, where a first-order and a second-order Hughes-like model are investigated. In particular, the authors find that the first-order model is not able to reproduce Braess' paradox, arguing that it is basically insensitive to (small) obstacles. Our results are not in contradiction with these ones since we reproduce the paradox using the basic behavior, whereas the authors of~\cite{twarogowska2014AMM} consider a model for rational pedestrians.

After this introduction, we can describe the organization of the paper: Section~\ref{sec:models} presents our new class of models; Section~\ref{sec:environ_optim} deals with the environment optimization; Section~\ref{sec:numerical_tests} presents exploratory but realistic numerical tests; finally, Section~\ref{sec:openpb} discusses analytical issues and sketches possible perspectives of further research.

\section{The models}
\label{sec:models}
The models we propose for pedestrian flows are first-order, two-dimensional, macroscopic, and composed by three interconnected elements: 
\paragraph{Small scale interactions} They describe short and medium-range influence of nearby walkers on a generic representative individual. They take into account the visual field of pedestrians and their sensory region~\cite{cristiani2011MMS,cristiani2014book}, namely a subset of the visual field where the presence of other people actually affects the walking dynamics.
\paragraph{Large scale interactions and path planning} They describe how pedestrians choose a path toward their destination according to their assumed predictive abilities (i.e., level of rationality).
\paragraph{Pedestrian-structure interactions} They describe the influence of the environment on pedestrian movements. This includes structural elements, such as obstacles, access or exit points.

The domain where pedestrians move, also called the \emph{walking area}, is a bounded set $\Omega\subset\R^2$. Possible obstacles are understood as holes in it, thus $\Omega$ identifies directly the portion of plane actually accessible to pedestrians.

We denote the time variable by $t$, the space variable by $x$, and the spatial density of the continuum representing the crowd by $\rho=\rho(t,x)$. The total number of pedestrians populating the walking area at time $t$ is $N_P(t)=\int_{\Omega}\rho(t,x)\,dx$.

Crowd dynamics in $\Omega$ obey the mass conservation principle, which implies that the density $\rho$ satisfies
\begin{equation}
	\partial_t\rho+\Div(\rho v)=0, \qquad\qquad \rho(0,\cdot)=\rho_0,
	\label{eq:conslaw}
\end{equation}
for any $t>0$ and $x\in\Omega$, $v=v(t,x)\in\R^2$ being the macroscopic velocity of the crowd and $\rho_0$ its initial distribution. Since, as stated in the Introduction, pedestrians are able to decide how to move by elaborating actively a behavioral strategy, we refrain from invoking inertial principles for obtaining their velocity $v$ (which would lead to a second-order model). We prefer to adopt a first-order modeling approach, in which pedestrian velocity is genuinely modeled by accounting for the superposition of two basic contributions: On the one hand, the desire to reach a specific destination, such as e.g., an exit point, which generates a \emph{behavioral velocity} $v_b^\ast$; On the other hand, the necessity to avoid crowded regions and collisions with other surrounding walkers, which generates an \emph{interaction velocity} $v_i$. Hence the total velocity is written as
\begin{equation}
	v=v_b^\ast+v_i.
	\label{eq:v}
\end{equation}

The interaction velocity $v_i$ is the component that pedestrians cannot control because it has to do with the effect produced by other walkers. In particular, it does not stem from choices made by individuals but rather from external conditions the latter undergo. The construction process will be studied in the next Section~\ref{sect:vi}.

Instead, the behavioral velocity $v_b^\ast$ is the component that pedestrians can control more directly, i.e., the one which they can effectively act on to implement their behavioral strategy. We have denoted it by an asterisk so as to convey the idea that pedestrians try to choose it in a way that they perceive as \emph{optimal} with respect to some personal preference criterion. The construction process of $v_b^\ast$ will be studied in detail in Section~\ref{sect:behavioral}. 

After being assembled as in~\eqref{eq:v}, the total velocity $v$ is finally corrected in order to take into account constraints imposed by the environment, such as obstacle walls or exit points, see Section~\ref{sect:obstacles}. 

\begin{remark}[Turning off the rationality]
\label{rem:turnoff}
Our model exhibits a clear separation between the controlled and the uncontrolled part of the velocity field~\eqref{eq:v}, making it easy to turn on or off the rationality. This feature is also present in the microscopic model for highly rational pedestrians~\cite{hoogendoorn2003OCAM} but not in other well-known macroscopic models, such as e.g.,~\cite{hoogendoorn2004TRBb,hughes2002TRB,xia2009PRE}.
\end{remark}

\subsection{Small scale interactions: Construction of $\boldsymbol{v_i}$}
\label{sect:vi}
\begin{figure}[!t]
\begin{center}
\includegraphics[width=0.2\textwidth]{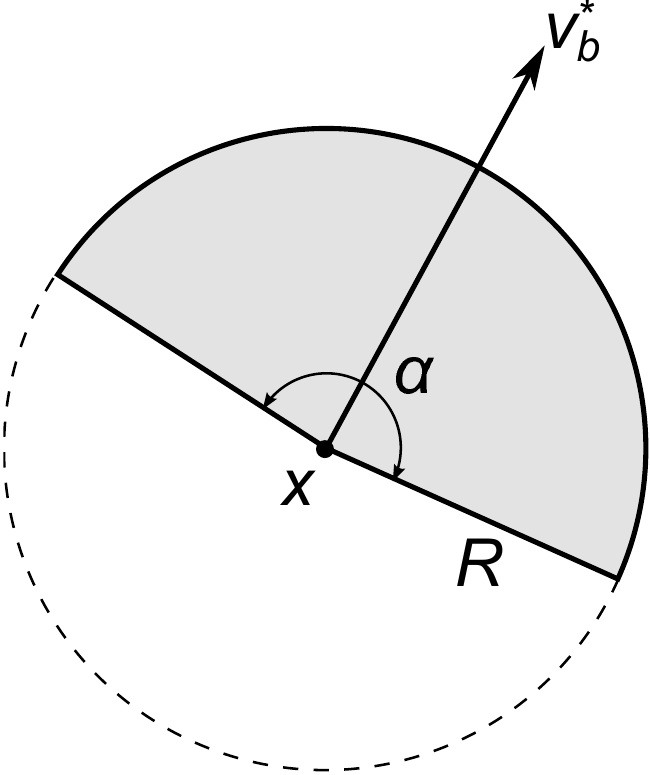}
\end{center}
\caption{The sensory region in the point $x$ oriented according to the behavioral velocity $v_b^\ast$.}
\label{fig:sensory_region}
\end{figure}

Following~\cite{cristiani2011MMS}, we express $v_i$ in nonlocal form:
\begin{equation}
	v_i=v_i[\rho(t,\cdot)](x)=\int_{\mathcal{S}(x)\cap\Omega}\mathcal{F}(y-x)\rho(t,y)\,dy
	\label{eq:vi}
\end{equation}
in order to model the ability of pedestrians to scan their surroundings and react to the distribution of other people therein. In particular, in~\eqref{eq:vi}:

\noindent $\bullet$ $\mathcal{S}(x)\subset\R^2$ is the \emph{sensory region} of the walker in $x$, i.e., the surrounding area within which s/he is sensitive to the presence of other pedestrians. Usually $\mathcal{S}(x)$ is assimilated to a circular sector of angular width $\alpha\sim 170^\circ$ -- the \emph{visual angle}, cf., e.g.,~\cite{robin2009TRB} -- while its radius $R$ -- the \emph{sensory radius} -- can vary from a few to several meters depending on the characteristic size of the environment under consideration. The orientation of $\mathcal{S}(x)$ is usually defined by the instantaneous direction of the vector $v_b^\ast$ in the point $x$, as illustrated in Fig.~\ref{fig:sensory_region}. This models the anisotropy of the interactions among pedestrians, who are sensitive to what is ``ahead'' with respect to their target direction of movement. In formulas:
$$ \mathcal{S}(x)=\left\{y\in\R^2\,:\,\abs{y-x}\leq R,\
	(y-x)\cdot v_b^\ast(x)\geq\cos{\frac{\alpha}{2}}\abs{y-x}\abs{v_b^\ast(x)}\right\}, $$
where $\cdot$ denotes the usual Euclidean inner product in $\R^2$;

\noindent $\bullet$ $\mathcal{F}:\R^2\to\R^2$ models the reaction of an individual in $x$ to another individual in $y$. For collision avoidance purposes, such a reaction is a repulsion, which we model as inversely proportional to the distance between the two individuals along the direction connecting them. Defining $r:=y-x$, we set $\mathcal{F}(r)=-Fr/\abs{r}^2$,
where $F>0$ is a constant tuning the strength of the repulsion. 
As suggested in~\cite{cristiani2012CDC}, a cut-off can be applied to $\mathcal{F}$ in order to avoid the singularity at $r=0$. Moreover, $v_i$ can be further regularized near $\partial\mathcal{S}(x)$ by introducing a suitable regularization of the characteristic function of $\mathcal{S}(x)$ in the integral~\eqref{eq:vi}.

\begin{remark}[Fundamental-diagram-free model]
\label{rem:fund_diag}
We stress that the pedestrian velocity~\eqref{eq:vi} is \emph{not} expressed as a function of the density by means of a fundamental diagram, as it often happens in the literature, see e.g.,~\cite{hoogendoorn2004TRBb,hughes2002TRB}. In our opinion this is a merit of our model, considering that fundamental diagrams are clearly a legacy of one-dimensional vehicular traffic models with no rigorous justification, from a kinematic point of view, in the case of two-dimensional flows. Conversely, by modeling the velocity field as in~\eqref{eq:v} with the nonlocal interaction term~\eqref{eq:vi} we are introducing a genuinely multidimensional model, which seems more natural for pedestrian flows.
\end{remark}

\subsection{Large scale interactions and path planning: Construction of $\boldsymbol{v_b^\ast}$}
\label{sect:behavioral}
In order to model pedestrian rationality as described in the Introduction, we need a mathematical tool able to capture the differences between basic, rational, and highly rational behaviors. To this end we consider the behavioral velocity field as a control function that each pedestrian chooses in a suitable set ${\cal A}$.

In more detail, we fix a performance criterion $\mathcal{C}$ and assume that pedestrians choose $v_b^\ast \in {\cal A}$ so as to get the best outcome from the path they will follow, i.e. $v_b^\ast$ is the $\mathcal{C}$-\emph{optimal} strategy.  Without loss of generality, in the following we assume that the maximal speed is normalized to $1$, so that ${\cal A}$ can be defined as the set of measurable control functions taking values in the closed unit ball $\overline{B_1(0)}$ (state constraints will be considered afterwards).

Given the crowd distribution $\rho$, individual pedestrians move according to the ordinary differential equation
\begin{equation}
	\dot{x}(t)=v_b^\ast(t,x(t))+v_i[\rho(t,\cdot)](x(t)),
	\label{eq:dyn_real}
\end{equation}
but in practice they might not have a clear picture of such dynamics. Thus, \emph{during the process of choice of their control} $v_b^\ast$ they might well be convinced to follow another dynamics, in general different from (\ref{eq:dyn_real}), given by
\begin{equation}
	\dot{y}(t)=f(t,y,v_b),\qquad v_b(\cdot)\in\mathcal A,
	\label{eq:dyn_perc}
\end{equation}
which in the following will be referred to as ``perceived dynamics'', due to their limited ability to observe the surroundings and/or predict the behavior of other walkers in $\Omega$. The precise expression of the ``perceived vector field'' $f$ is exactly what characterizes the different degrees of rationality and will be discussed in detail in the following sections. Here, it is important to stress that pedestrians do choose the field $v_b^\ast$ which is $\mathcal{C}$-optimal along the trajectories of~\eqref{eq:dyn_perc} but which could be not as optimal when the actual evolution~\eqref{eq:dyn_real} takes place.

The choice of the performance criterion $\mathcal{C}$ depends on the situation to be modeled. To fix the ideas, we focus our attention on two different frameworks: Evacuation problems and fixed horizon evolutions. 

In \emph{evacuation problems} we assume that pedestrians have a target $\Sigma$ (= exits from $\Omega$) on the boundary $\partial \Omega$ of the domain, and that the cost criterion to be minimized is given by the time needed to reach $\Sigma$ from any initial position $x\in\Omega$. This corresponds to choosing\footnote{Here and henceforth solutions to~\eqref{eq:dyn_perc} in $[a,b]\subset\R$ are understood in Carath\'{e}odory sense.
}:
\begin{equation}
\begin{split}
\mathcal{C}(x,v_b(\cdot)):=\min\big\{\mathfrak{t}\geq 0\,:\,\exists\,y(\cdot)\ \text{solution to } \dot{y}(t)=f(t,y,v_b) \text{ in\ } 
	[0,\mathfrak{t}]\  \text{s.t.\ }\\ y(0)=x,\,y(\mathfrak{t})\in\Sigma\big\}, \quad x\in\Omega, \quad v_b(\cdot)\in\mathcal A.
\label{eq:min_time}
\end{split}
\end{equation}

However, as pointed out in~\cite{mason2012SICON}, evacuation problems ruled by the performance criterion~\eqref{eq:min_time} are in some cases not enough to describe crowd behaviors observed in real situations. Even if pedestrians have targets, and even if the density remains sufficiently low in the whole domain so as to avoid deviations due to crowded regions, the chosen paths are often not time-minimal. This suggests that the subjective performance criterion actually adopted by pedestrians may be somehow different from~\eqref{eq:min_time} and motivates us to consider also fixed horizon evolutions.

In \emph{fixed horizon evolutions} pedestrian dynamics take place over a fixed time interval $[0,T]$ and the performance criterion consists in a running cost, accounting for the cost of a specific path, and an exit cost, accounting for the possible failure to reach the target at time $T$. Namely, we consider
\begin{equation}
	\mathcal{C}(t,x,v_b(\cdot)):=\int_t^T\ell\big(y(\sigma),v_b(\sigma)\big)\,d\sigma+g(y(T)), \quad (t,x)\in [0,T]\times\Omega,\quad v_b(\cdot)\in\mathcal A
	\label{eq:fixed_horiz}
\end{equation}
where $\ell$ is a positive, continuous, and bounded running cost, $g$ is a continuous and bounded below exit cost, and $y(\cdot)$ is again a solution to $\dot{y}(\sigma)=f(\sigma,y,v_b)$ in $[t,T]$ with initial datum $y(t)=x$. Notice that this formulation allows one to model also unhurried pedestrians, for instance when their travel purpose is leisure or when their goal is an appointment at a later time $T$.

For both frameworks~\eqref{eq:min_time},~\eqref{eq:fixed_horiz} we now describe how to choose the perceived vector field $f$ in~\eqref{eq:dyn_perc} so as to model the basic, rational, and highly rational behaviors. Moreover, we introduce a new family of hybrid models, still described in terms of suitable perceived vector fields, where the rationality degree can be freely tuned.

\subsubsection{Basic behavior}
In case of \emph{basic behavior} we simply assume that each pedestrian chooses $v_b^\ast$ without considering what other walkers do in the meanwhile. Hence an individual in $x\in\Omega$ at time $t$ chooses his/her strategy $v_b^\ast$ only on the basis of his/her knowledge of the walking area and of the chosen cost $\mathcal{C}$, ignoring completely the crowd distribution $\rho$. Formally, this implies that the control value $v_b^\ast$ is chosen so as to be $\mathcal{C}$-optimal when the perceived dynamics~\eqref{eq:dyn_perc} are given by 
$\dot{y}(t)=v_b(t)$,
with $v_b(\cdot)\in\mathcal A$ and no interaction contribution.

When dealing with evacuation problems, the construction of $v_b^\ast$ can proceed as follows. We define the \emph{value function} of the problem as
\begin{equation}
	\phi(x):=\min_{v_b(\cdot)\in\mathcal{A}}\mathcal{C}(x,v_b(\cdot)), \qquad x\in\Omega,
	\label{eq:basic_vf}
\end{equation}
where $\mathcal{C}$ was given in~\eqref{eq:min_time} with $f(t,y,v_b)=v_b$, and we assume for simplicity that $\phi(x)<+\infty$ for all $x\in\Omega$. Notice that $\phi(x)$ is nothing but the minimum time needed to reach $\Sigma$ from any initial position $x\in\Omega$. Then, it is well known in control theory that $\phi$ is a bounded uniformly continuous viscosity solution to the eikonal equation
\begin{equation}
	\abs{\nabla\phi(x)}-1 = 0, \quad x\in\Omega
	\qquad\text{with}\qquad
	\phi(x)=0, \quad x\in\Sigma
	\label{eq:eikonal_equation}
\end{equation}
(we refer to~\cite[Sect.\ 8.4.5]{falconebook} for the correct treatment of the conditions on the boundary $\partial\Omega\setminus\Sigma$). After finding $\phi$ from~\eqref{eq:eikonal_equation}, the optimal behavioral velocity field can be taken in feedback form as 
\begin{equation}
	v_b^\ast(x)=-\frac{\nabla\phi(x)}{\abs{\nabla\phi(x)}}, \qquad x\in\Omega.
	\label{eq:basic_optimal1}
\end{equation}
Notice that $v_b^\ast$ is time-independent coherently with the fact that pedestrians are insensitive to the evolution of the crowd.

When dealing with fixed horizon evolutions, the procedure is similar. We define the value function for the problem as
\begin{equation}
	\phi(t,x):=\min_{v_b(\cdot)\in\mathcal{A}}\mathcal{C}(t,x,v_b(\cdot)), \qquad (t,x)\in [0,T]\times\Omega,
	\label{eq:fixed_horiz_vf}
\end{equation}
where $\mathcal{C}$ was given in~\eqref{eq:fixed_horiz} with $f(t,y,v_b)=v_b$, and we rely on control theory to find $\phi$ as the unique viscosity solution of the following Hamilton-Jacobi-Bellman (HJB) equation on $(0,T)\times\Omega$:
\begin{equation}
	-\partial_t\phi(t,x)+\max_{v_b\in\overline{B_1(0)}}\left\{-v_b\cdot\nabla{\phi(t,x)}-\ell(x,v_b)\right\}=0
	\label{eq:basic_hjb}
\end{equation}
with final condition $\phi(T,x)=g(x)$ and suitable boundary conditions (see again~\cite{falconebook}). After finding $\phi$ from~\eqref{eq:basic_hjb}, the optimal behavioral velocity $v_b^\ast$ can be selected invoking Pontryagin maximum principle (see~\cite{bardibook}) by constructing a feedback control such that
\begin{equation}
	v_b^\ast(t,x)\in\argmax_{v_b\in \overline{B_1(0)}}\left\{-v_b\cdot\nabla{\phi(t,x)}-\ell(x,v_b)\right\}
		\qquad \forall\,(t,x)\in [0,T]\times\Omega.
	\label{eq:basic_optimal2}
\end{equation}


\begin{remark}
Observe that, even if the vector field $v_b^\ast$ is computed independently of time or if it is assigned \emph{a priori}, crowd dynamics do not reduce to plain movements along the trails given by the integral curves of $v_b^\ast$. Indeed the interaction term $v_i[\rho]$ (cf. the previous Section~\ref{sect:vi}) is still present in~\eqref{eq:v} and affects both the real dynamics~\eqref{eq:dyn_real} and the continuity equation~\eqref{eq:conslaw}. This is an advantageous effect of the clear separation between the controlled part (influenced by rationality) and the uncontrolled part of the velocity field~\eqref{eq:v}, already pointed out in Remark~\ref{rem:turnoff}.
\end{remark}

\subsubsection{Rational behavior}
In case of \emph{rational behavior} we assume that at a fixed time $\tau\geq 0$ each pedestrian is aware of the distribution $\rho(\tau,\cdot)$ in $\Omega$ and wants to use this information to choose the optimal path to his/her destination. Then, treating $\tau$ as a fixed parameter, we assume that the perceived dynamics~\eqref{eq:dyn_perc} for rational pedestrian are given by 
$$ \dot{y}(t)=v_b(t)+v_i[\rho(\tau,\cdot)](y(t)), \qquad v_b(\cdot)\in\mathcal A.$$ 
The nonlocal vector field $v_i$, defined in~\eqref{eq:vi}, does not depend on time in this case, because the density acknowledged by pedestrians is ``frozen'' at time $\tau$. Hence the individuals can adopt the same construction as before for the feedback strategy $x\mapsto v_b^\ast(\tau,x)$ in $\Omega$, with minor adaptations. Finally, we require rational pedestrians to solve~\eqref{eq:conslaw}--\eqref{eq:v} by \emph{repeating at each time} $t=\tau$ such a construction of $v_b^\ast(\tau,\cdot)$.

For completeness, we now sketch the precise construction of $\mathcal{C}$-optimal strategies $v_b^\ast$ in the two scenarios we are interested in. With respect to the basic behavior, the main difference is that now the criterion $\mathcal{C}$ depends also on the fixed parameter $\tau\geq 0$.

In the case of evacuation problems, the cost becomes
\begin{equation*}
\begin{split}
\mathcal{C}\big(x,v_b(\cdot);\tau,\rho(\tau,\cdot)\big)\!:=\!\min\!\big\{\mathfrak{t}\geq\tau\,:\,\exists\,y(\cdot)\,\text{solution to }\dot{y}(t)=v_b+v_i[\rho(\tau,\cdot)](y) \text{ in }\, \\
[\tau,\mathfrak{t}]\,\text{ s.t. } y(\tau)=x,\,y(\mathfrak{t})\in\Sigma\big\}, \quad x\in\Omega,\quad v_b(\cdot)\in\mathcal A.
\end{split}
\end{equation*}
Defining the value function $\phi(x)=\phi_\tau(x)$ as in~\eqref{eq:basic_vf}, and still assuming that $\phi_\tau(x)<+\infty$ in $\Omega$, we can find $\phi$ as a bounded uniformly continuous viscosity solution to:
\begin{equation}
	\abs{\nabla{\phi_\tau(x)}}-v_i[\rho(\tau,\cdot)](x)\cdot\nabla{\phi_\tau(x)}-1=0, \qquad x\in\Omega
	\label{eq:HJB_rational1}
\end{equation}
with boundary conditions $\phi_\tau=0$ on $\Sigma$. After finding $\phi_\tau$ from~\eqref{eq:HJB_rational1}, the $\mathcal{C}$-optimal feedback can be chosen as
$$ v_b^\ast(\tau,x)=-\frac{\nabla{\phi_\tau(x)}}{\abs{\nabla{\phi_\tau(x)}}}, \qquad x\in\Omega, $$
which is completely analogous to~\eqref{eq:basic_optimal1}. Notice that, at fixed $\tau\geq 0$, the HJB equation~\eqref{eq:HJB_rational1} is independent of the continuity equation~\eqref{eq:conslaw}--\eqref{eq:v} because the vector field $x\mapsto v_i[\rho(\tau,\cdot)](x)$, thus also the optimization procedure used to compute $v_b^\ast$, does not require the knowledge of $\rho$ at any time $t\neq \tau$.

In the case of fixed horizon problems, the situation is completely analogous. We treat $\tau$ as a parameter and define, for $(t,x)\in [\tau,T]\times\Omega$, both a cost functional analogous to~\eqref{eq:fixed_horiz} and a value function $\phi(t,x)=\phi_\tau(t,x)$ analogous to~\eqref{eq:fixed_horiz_vf}. Then $\phi_\tau$ can be computed by solving in $(\tau,T)\times\Omega$ the HJB equation
\begin{equation}
	-\partial_t\phi_\tau(t,x)+\max_{v_b\in\overline{B_1(0)}}
		\left\{-\left(v_b+v_i[\rho(\tau,\cdot)](x)\right)\cdot\nabla{\phi_\tau(t,x)}-\ell(x,v_b)\right\}=0,
	\label{eq:HJB_rational2}
\end{equation}
with terminal condition $\phi_\tau(T,x)=g(x)$. Once again, this equation is independent of the continuity equation. Finally, $v_b^\ast(\tau,\cdot)$ is selected, like in~\eqref{eq:basic_optimal2}, among the values that realize the maximum of the Hamiltonian in~\eqref{eq:HJB_rational2} at time $t=\tau$.

\subsubsection{Highly rational behavior}
The procedure described in the previous section for rational pedestrians basically states that individuals, when choosing $v_b^\ast$ at a certain time $t$, merely react to the crowd distribution at that precise moment. Therefore, even if they continuously adapt their strategies in time, the resulting behavior is not globally optimal. In other words, rational pedestrians make choices that can be good instantaneously but not really rewarding in the long run.

In order to introduce predictive abilities in our model, we have to choose appropriately the perceived vector field~\eqref{eq:dyn_perc}. Different choices are required in the two scenarios we focus on. For fixed horizon problems, we assume that pedestrians are aware of the crowd distribution $\rho$ in the whole time-space $[0,T]\times\Omega$, namely that they seek a global optimum in both space and time. This will lead us naturally to Nash equilibria and mean field games. Conversely, for evacuation problems we need to extend the minimum time problem, which is intrinsically autonomous, to time-space settings, see, e.g.,~\cite{bokanowski2011IFAC}. Let us discuss separately the details for the two cases.

\paragraph{Fixed horizon problems}
In this case we assume that pedestrians are able to predict exactly the evolution of the density $\rho$ during the whole time horizon $[0,T]$ where crowd dynamics take place. In other words, we assume that they can assess the long term effect of their choices on the other walkers and that they choose accordingly the strategy $(t,x)\in[0,T]\times\Omega\mapsto v_b^\ast(t,x)$ so as to get their best outcome at time $T$. With this behavior we are thus describing perceived dynamics~\eqref{eq:dyn_perc} of the form
$$ \dot{y}(t)=v_b(t)+v_i[\rho(t,\cdot)](y(t)),\qquad v_b(\cdot)\in\mathcal A, $$
where the interaction velocity field has now to be regarded as a function of time and space because $\rho$ is not ``frozen'' as it is in the rational behavior. Note that, in this case, once the optimal behavioral velocity is taken the perceived and the real dynamics coincide. As a cost functional we consider $\mathcal{C}(t,x,v_b(\cdot))=\mathcal{C}(t,x,v_b(\cdot)\,;\rho(\cdot,\cdot))$ of the form~\eqref{eq:fixed_horiz} and assume that pedestrians construct the $\mathcal{C}$-optimal feedback $v_b^\ast$ along the lines of the previous cases. They solve an HJB equation for the value function~\eqref{eq:basic_hjb}, which takes the form:
\begin{equation}
	-\partial_t\phi(t,x)+\max_{v_b\in\overline{B_1(0)}}\left\{-(v_b+v_i[\rho(t,\cdot)](x))\cdot\nabla{\phi(t,x)}
		-\ell(x,v_b)\right\}=0,
	\label{eq:hrational_HJB_FH}
\end{equation}
equipped with the final condition $\phi(T,x)=g(x)$. After finding $\phi$ from~\eqref{eq:hrational_HJB_FH}, the optimal feedback control $v_b^\ast(t,x)$ is recovered like in~\eqref{eq:basic_optimal2}. Notice that this time the HJB equation is coupled to the continuity equation~\eqref{eq:conslaw}--\eqref{eq:v}, because the optimization procedure used to compute $v_b^\ast$ requires the knowledge of $\rho$ at all times. Therefore modeling the highly rational behavior implies dealing with the following \emph{forward-backward} system of PDEs:
\begin{equation}
\begin{cases}
	\partial_t\rho(t,x)+\Div\!\Big(\big(v_b^\ast(t,x)+v_i[\rho(t,\cdot)](x)\big)\rho(t,x)\Big)=0 \\[3mm]
	-\partial_t\phi(t,x)+
		\displaystyle{\max_{v_b\in\overline{B_1(0)}}}\left\{-(v_b+v_i[\rho(t,\cdot)](x))\cdot\nabla{\phi(t,x)}-\ell(x,v_b)\right\}=0 \\[3mm]
	v_b^\ast(t,x)\in
		\displaystyle{\argmax_{v_b\in\overline{B_1(0)}}}\left\{-(v_b+v_i[\rho(t,\cdot)](x))\cdot\nabla{\phi(t,x)}-\ell(x,v_b)\right\}
\end{cases}
\label{eq:mfg}
\end{equation}
with initial condition $\rho(0,\cdot)=\rho_0$ and final condition $\phi(T,\cdot)=g$ in $\Omega$. The system is \emph{fully coupled} because of the presence of the whole density $\rho$ in the second equation as well as of $v_b^\ast$ in the first one.

Mathematically, such a system can be interpreted as a first-order mean field game along the lines of~\cite{LL}, see also~\cite{cardanotes}. This allows us to justify rigorously the choice of the terminology ``highly rational'' for describing such a behavior. Indeed, a solution pair $(\rho,\phi)$, when it exists, represents the crowd distribution and the corresponding performance outcome for pedestrians who reach a Nash equilibrium among themselves. 

\paragraph{Evacuation problems}
In this case we need to recast the problem in a time-space domain, so as to be able to include the evolution of the density $\rho$ into the problem, with minor modifications to the equations themselves.
To this end, we define the \emph{time-space domain} and the \emph{time-space target} by setting
$\ext{\Omega}:=\{t>0\}\times\Omega$ and $\ext{\Sigma}:=\{t\geq 0\}\times\Sigma$.
We also extend the velocity field as
$$ \ext{v}(t,x):=
	\begin{pmatrix}
		1 \\
		v(t,x)
	\end{pmatrix}
	=
	\begin{pmatrix}
		1 \\
		v_b(t)+v_i[\rho(t,\cdot)](x)
	\end{pmatrix}, 
	\qquad (t,x)\in\ext{\Omega},\qquad v_b(\cdot)\in\mathcal A,
$$
and we assume that the perceived dynamics for the new variable $z\in\ext{\Omega}$ are given by
\begin{equation}
	\dot{z}(s)=\ext{v}(z(s)).
	\label{eq:dyn_extended}
\end{equation}
In other words, the time variable is seen as an additional space variable and a fictitious time $s$ is used to move pedestrians in the time-space. Note that pedestrians move in time with fixed speed $1$, namely they are passively carried toward the future (time traveling will be attacked in a forthcoming paper). To complete the description of the behavior, we assume that highly rational pedestrians choose $v_b^\ast$ by solving an ordinary minimum time problem on $\ext{\Omega}$ with target $\ext{\Sigma}$. Since advection in time is not controlled, any optimal control for this extended problem corresponds to a time-optimal control in $\Omega$ with target $\Sigma$.

With respect to the case of rational pedestrians, the main difference is that here the whole evolution of $\rho$ in $\ext{\Omega}$ is needed for the construction of $v_b^\ast$. Let us define
\begin{equation*}
\begin{split}
\mathcal{C}(t,x,v_b(\cdot);\rho(\cdot,\cdot))\!:=\!
	\min\big\{\mathfrak{s}\geq 0\,:\,\exists\,z(\cdot)\,\text{solution to~\eqref{eq:dyn_extended} in } [0,\mathfrak{s}] \text{ s.t.} \\ z(0)=(t,x), z(\mathfrak{s})\in\ext{\Sigma}\big\},\quad (t,x)\in\ext{\Omega},\quad v_b(\cdot)\in\mathcal A
\end{split}
\end{equation*}
and the corresponding value function like in~\eqref{eq:fixed_horiz_vf}. The associated HJB equation is given by
\begin{equation}
	\max_{v_b\in\overline{B_1(0)}}\left\{-\left(1,\,v_b+v_i[\rho(t,\cdot)](x)\right)\cdot
		\nabla_{t,x}{\phi(t,x)}-1\right\}=0, \quad (t,x)\in\ext{\Omega}\setminus\ext{\Sigma},
	\label{eq:HJB_hrational_MT}
\end{equation}
with boundary conditions $\phi(t,x)=0$ on $\ext{\Sigma}$, $\nabla_{t,x}$ denoting the gradient w.r.t.\ both time and space variables. After finding $\phi$ from~\eqref{eq:HJB_hrational_MT}, the optimal feedback control $v_b^\ast(t,x)$ can be chosen as usual among the values which realize the maximum in~\eqref{eq:HJB_hrational_MT}. Consequently, modeling the highly rational behavior requires solving the following \emph{forward-backward} system of PDEs:
\begin{equation}
	\begin{cases}
		\partial_t \rho(t,x)+\Div\!\Big(\big(v_b^\ast(t,x)+v_i[\rho(t,\cdot)](x)\big)\rho(t,x)\Big)=0 \\[3mm]
		\displaystyle{\max_{v_b\in\overline{B_1(0)}}}\left\{-\left(1,\,v_b+v_i[\rho(t,\cdot)](x)\right)\cdot\nabla_{t,x}{\phi(t,x)}-1\right\}=0 \\[3mm]
		v_b^\ast(t,x)\in\displaystyle{\argmax_{v_b\in\overline{B_1(0)}}}\left\{-\left(1,\,v_b+v_i[\rho(t,\cdot)](x)\right)\cdot
			\nabla_{t,x}{\phi(t,x)}-1\right\}
	\end{cases}
	\label{eq:hrational_mt}
\end{equation}
for $t>0$, $x\in\Omega$, with initial condition $\rho(0,\cdot)=\rho_0$ in $\Omega$ and boundary condition $\phi(t,x)=0$ on $\ext{\Sigma}$. The second equation in~\eqref{eq:hrational_mt} is said to be \emph{backward} because characteristic curves flow backward in time, from the future to the past. Note that, like in the case of system~\eqref{eq:mfg}, system~\eqref{eq:hrational_mt} is fully coupled.

\subsubsection{$\theta$-rational behavior}
Using the machinery introduced above, we can also define hybrid behaviors between the rational and the highly rational ones in both scenarios we are interested in. To this end, let us introduce a new parameter $\theta$, which is taken in $[0,T]$ for fixed horizon problems and in $[0,+\infty)$ for evacuation problems. Such a parameter serves to span the whole family of behaviors, from the rational one, corresponding to $\theta=0$, to the highly rational one, corresponding to $\theta=T$ (resp., $\theta\to +\infty$) if the finite horizon cost~\eqref{eq:fixed_horiz} (resp., the minimum time cost~\eqref{eq:min_time}) is employed.

In more detail, we define the \textit{$\theta$-rational behavior} by assuming that pedestrians do have predictive abilities however limited in time, in particular extending only up to a time $\theta$ in the future. To set the model, we adopt a strategy similar to the one we used for rational pedestrians. Namely, we prescribe a procedure for constructing the feedback $v_b^\ast(\tau,\cdot)$ in $\Omega$ at any fixed time $\tau\geq 0$, then for every time $t$ we repeat the process with $\tau=t$. 

Thus, fixing any time $\tau$, we do not only assume that pedestrians are aware of the current density distribution $\rho(\tau,\cdot)$ in $\Omega$ but also that they can forecast the evolution of $\rho$ up to time $\tau+\theta$. In other words, we assume that pedestrians are aware of the distribution $\rho^\theta_\tau:[\tau,\tau+\theta]\times\Omega\to\R$ resulting from~\eqref{eq:conslaw} with initial datum $\rho(\tau,\cdot)$ at time $t=\tau$ and that they prolong it by setting $\rho^\theta_\tau(t,x):=\rho^\theta_\tau(\tau+\theta,x)$ for every $t>\tau+\theta$ and $x\in\Omega$. Using such an auxiliary density $\rho^\theta_\tau\colon[\tau,+\infty)\times\Omega\to\R$, pedestrians adopt the following perceived dynamics~\eqref{eq:dyn_perc}:
$$ \dot{y}(t)=v_b(t)+v_i[\rho^\theta_{\tau}(t,\cdot)](y(t)), \qquad t>\tau,\qquad v_b(\cdot)\in\mathcal A,$$
which they use to choose the behavioral field $v_b^{\ast,\theta}$ entering in turn the equation for $\rho^\theta_\tau$. The model construction is analogous to that performed for highly rational pedestrians and finally leads to a \emph{fully coupled} system very similar to~\eqref{eq:mfg} (resp.,~\eqref{eq:hrational_mt}) but defined on the time interval $[\tau,T]$ (resp., $[\tau,+\infty)$) instead of $[0,T]$ (resp., $[0,+\infty)$) and featuring a continuity equation for the auxiliary density $\rho^\theta_\tau$ rather than for $\rho$. For brevity we provide details only for the evacuation problem, the fixed horizon problem being completely analogous. In this case, by defining $\ext{\Omega}_\tau:=\{t>\tau\}\times\Omega$ and $\ext{\Sigma}_\tau:=\{t\geq\tau\}\times\Sigma$, we obtain the following system of PDEs:
\begin{equation}
	\begin{cases}
		\partial_t\rho^\theta_\tau+\Div\big(\rho^\theta_\tau(v_b^{\ast,\theta}+v_i[\rho^\theta_\tau])\big)=0 & \text{in\ } (\tau,\tau+\theta)\times\Omega \\[3mm]
		\displaystyle{\max_{v_b\in\overline{B_1(0)}}}\left\{-\left(1,\,v_b+v_i[\rho^\theta_\tau(t,\cdot)](x)\right)\cdot
			\nabla_{t,x}{\phi^\theta_\tau(t,x)}-1\right\}=0 & \text{in\ } \ext{\Omega}_\tau \\[3mm]
		v_b^{\ast,\theta}(t,x)\in\displaystyle{\argmax_{v_b\in\overline{B_1(0)}}}\left\{-\left(1,\,v_b+v_i[\rho^\theta_\tau(t,\cdot)](x)\right)\cdot
			\nabla_{t,x}{\phi^\theta_\tau(t,x)}-1\right\} & \text{in\ } \ext{\Omega}_\tau
	\end{cases}
	\label{eq:theta_rat}
\end{equation}
with additionally $\rho^\theta_\tau(t,x)=\rho^\theta_\tau(\tau+\theta,x)$ in $\{t\geq\tau+\theta\}\times\Omega$, to extend $\rho^\theta_\tau$ to the whole domain $\ext{\Omega}_\tau$, and conditions $\rho^\theta_\tau(\tau,x)=\rho(\tau,x)$ in $\Omega$ and $\phi^\theta_\tau(t,x)=0$ in $\ext{\Sigma}_\tau$ for the forward-backward equations.

Whenever a solution pair $(\rho^\theta_\tau,\phi^\theta_\tau)$ to~\eqref{eq:theta_rat} is found, the feedback $(t,x)\mapsto v_b^{\ast,\theta}(t,x)$ is defined for all $(t,x)\in[\tau,+\infty)\times\Omega$ and we can just set $v_b^\ast(\tau,x)=v_b^{\ast,\theta}(\tau,x)$ for $x\in\Omega$. The process is then repeated for every time $t$, with $\tau=t$. 

\subsection{Pedestrian-structure interactions}
\label{sect:obstacles}
The influence of structural elements, mainly obstacles, on pedestrian behavior is modeled in terms of boundary conditions for the total velocity $v=v_b^\ast+v_i$. We recall that obstacles are understood as holes of $\Omega$, therefore their edges are inner boundaries of the domain. Denoting by $n=n(x)$ an outer normal vector to $\partial\Omega$ in the point $x$, we impose the ``impermeability'' condition:
\begin{equation}
	v\cdot n\leq 0,
	\label{eq:bc_v}
\end{equation}
meaning that the total velocity given by~\eqref{eq:v} has to be corrected near an obstacle or near the outer boundary of $\Omega$ so that the crowd scrapes against the walls without penetrating them. In the case of a piecewise regular boundary $\partial \Omega$, we require~\eqref{eq:bc_v} to hold for all vectors $n$ in the (Bouligand) tangent cone to $\Omega$ in the point $x\in\partial\Omega$.

On the whole, what we do technically is to project the velocity field~\eqref{eq:v} onto the space of admissible fields, which are precisely those fulfilling~\eqref{eq:bc_v} on $\partial\Omega$. With this choice, the domain is \emph{viable}, i.e., any pedestrian trajectory starting from $\overline{\Omega}$ remains in $\overline{\Omega}$ for positive times and cannot cross any portion of the boundary $\partial\Omega$. 

\subsection{Nondimensionalization of the equations}
It is useful to rewrite model \eqref{eq:conslaw}--\eqref{eq:vi} in dimensionless form so as to work with equations which are independent of the orders of magnitude of the physical variables. Some of such variables indeed vary from case to case (for instance the size of the domain and the sensitivity radius), while dimensionless equations can be implemented and solved always in the same manner, converting back dimensionless to dimensional results \emph{a posteriori}.

Let $L$, $V$, $\varrho$ be characteristic values of length, speed, and density, whence we deduce also the characteristic time $L/V$. For instance, the engineering literature indicates that pedestrians in normal conditions (i.e., no panic) walk at $V=O(1\unit{m/s})$, which corresponds to a packing $\varrho=O(1\unit{ped/m^2})$. We introduce the following dimensionless variables, functions, and parameters:
$$ \adim{x}=\frac{x}{L}, \quad \adim{t}=\frac{V}{L}t, \quad \adim{v}(\adim{t},\adim{x})=\frac{1}{V}v\left(\frac{L}{V}\adim{t},L\adim{x}\right),
	\quad \adim{\rho}(\adim{t},\adim{x})=\frac{1}{\varrho}\rho\left(\frac{L}{V}\adim{t},L\adim{x}\right), \quad
	\adim{F}=\frac{F\varrho L}{V} $$
to discover that~\eqref{eq:conslaw} becomes
$\partial_{\adim{t}}\adim{\rho}+\adim{\Div}(\adim{\rho} \adim{v})=0$,
where $\adim{\Div}$ denotes the divergence operator with respect to $\adim{x}$. Moreover, from~\eqref{eq:v}--\eqref{eq:vi} we get
$$ \adim{v}(\adim{t},\adim{x})=\adim{v}_b^\ast(\adim{t},\adim{x})
	+\int_{\adim{\mathcal{S}}(\adim{x})\cap\adim{\Omega}}\adim{\mathcal{F}}(\adim{y}-\adim{x})\adim{\rho}(\adim{t},\adim{y})\,d\adim{y}, 
	$$
where $\adim{\mathcal{F}}(\adim{r})=-\adim{F}\adim{r}/\abs{\adim{r}}^2$, and $\adim{\mathcal{S}}(\adim{x})=\frac{1}{L}\mathcal{S}(L\adim{x})$ and $\adim{\Omega}=\frac{1}{L}\Omega$ are the dimensionless $L$-homothetic versions of the sensory region and the walking area, respectively. 

In the rest of the paper we will always consider the equations in their dimensionless form, but to avoid unnecessarily heavy notations we will drop the symbol $\tilde{}$ over the dimensionless variables.

\subsection{Numerical approximation}
The conservation law~\eqref{eq:conslaw} is discretized by means of the scheme firstly proposed in~\cite{piccoli2011ARMA} and then used extensively in~\cite{cristiani2011MMS,cristiani2012CDC,cristiani2014book}. It is a truly two-dimensional first-order reasonably fast conservative scheme, which has been proved to converge to a weak solution of~\eqref{eq:conslaw}, see~\cite{piccoli2013AAM,TosinFrasca}, and to describe adequately the main features of pedestrian flow, including merging and splitting, although it exhibits a non-negligible numerical diffusion (its one-dimensional version coincides with the classical upwind scheme). 

The HJB equations~\eqref{eq:eikonal_equation} and~\eqref{eq:HJB_rational1} are discretized by means of an iterative first-order semi-Lagrangian scheme. The interested reader can find a complete introduction to the topic in the recent book~\cite{falconebook} (see also~\cite{bardibook}). The scheme used here is described (with complete references) in~\cite{cacace2014SISC}. The unit ball $\overline{B_1(0)}$ is discretized with $32$ points, all placed on the boundary. The reconstruction of the values of the solution at non-mesh points is obtained by means of a three-point linear interpolation. The Fast Sweeping technique (see again~\cite{cacace2014SISC} for explanations and references) is used to speed up the convergence. The HJB equation~\eqref{eq:HJB_hrational_MT} is discretized by means of a generalization of the scheme used for~\eqref{eq:HJB_rational1}, extending both the grid and the interpolation to the time dimension, cf.\ \cite{carlini2014SINUM}.

For all equations the space domain is discretized by means of a structured grid with $100\times 100$ nodes. 

\section{Environment optimization}
\label{sec:environ_optim}
After modeling the various pedestrian behaviors, we can use them to control the movement of the crowd, steering the dynamics toward a desired behavior. To this end, we have to define the \emph{natural behavior} of pedestrians, i.e., the behavior they are expected to assume in reality, and a \emph{target behavior}, i.e., a desired behavior we would like them to assume, which is considered optimal with respect to some given criterion. The natural behavior will be typically the irrational/basic/rational one, while the target behavior will be typically the rational/highly rational/optimal one. Correspondingly, for $t>0$ and $x\in\Omega$, we denote by $\rho^{\textup{n}}(t,x)$ (resp., $\rho^{\textup{t}}(t,x)$) the density evolution according to the natural (resp., target) behavior.

\subsection{Main ingredients}\label{sec:control_and_cost} 
The two main ingredients of the environment optimization problem are the \emph{environmental control} and the \emph{environmental cost}. They are defined as follows.
\paragraph{The environmental control} 
We assume that one can introduce in the domain additional obstacles denoted by $\mathcal{O}_\lambda$, where $\lambda\in\R^{N_\mathcal{O}}$ for some $N_\mathcal{O}\in\mathbb N$ is any minimal set of parameters which determines univocally these obstacles. Thus, the spatial domain that pedestrians can actually access is $\Omega\setminus\mathcal{O}_\lambda$. Note that the set of parameters will be in general constrained, hence we will have $\lambda\in\Lambda\subset\R^{N_\mathcal{O}}$ for a suitable set $\Lambda$ of admissible parameters. Constraints can be dictated by either modeling needs (e.g., the fact that obstacles be completely contained in $\Omega$), compatibility conditions (such as, e.g., excluding obstacles overlapping with the initial crowd distribution $\rho_0$), or practical convenience (e.g., too large obstacles might not be desirable). In the following, we will say that an obstacle $\mathcal{O}_\lambda$ is \emph{admissible} if $\lambda\in\Lambda$.
\paragraph{The environmental cost} 
By varying $\lambda$ in $\Lambda$ we modify the environment in such a way that the natural behavior in the new controlled environment is as close as possible to the target one in the original environment. For any $t>0$ and $x\in\Omega\setminus\mathcal{O}_\lambda$, let $\rho^{\textup{n},\lambda}(t,x)$ be the density evolution according to the natural behavior in the controlled environment.  The cost criterion $\Delta(\lambda)$, used to measure the distance between the two behaviors, and to be minimized, can be defined in several ways. Here we list three of them.

\noindent 1. We define the evacuation time according to the controlled natural behavior as 
\begin{equation}
	\tevac^{\textup{n},\lambda}:=\min\{t>0\,:\,\rho^{\textup{n},\lambda}(t,x)=0\ \forall\,x\in\Omega\setminus\mathcal{O}_\lambda\},
	\label{def:tevac_natural}
\end{equation}
and the evacuation time according to the target behavior as 
\begin{equation}
	\tevac^{\textup{t}}:=\min\{t>0\,:\,\rho^{\textup{t}}(t,x)=0\ \forall\,x\in\Omega\}.
	\label{def:tevac_target}
\end{equation}
Then we define $\Delta_1$ to be the difference between the two evacuation times, i.e.,
\begin{equation}
	\Delta_1(\lambda):=\abs{\tevac^{\textup{n},\lambda}-\tevac^{\textup{t}}}.
	\label{def:Delta1}
\end{equation}

\noindent 2. Let $e_1,\,e_2,\,\dots,\,e_{N_e}$ be the exits of the domain, whose total number is $N_e$, and $P(e_k)$ the number of pedestrians leaving the domain through the exit $e_k$, for $k=1,\,\dots,\,N_e$. We introduce the vector $P_e:=\big(P(e_1),\,P(e_2),\,\dots,\,P(e_{N_E})\big)$ and define $\Delta_2$ as
\begin{equation}
	\Delta_2(\lambda):=\norm{(P_{e}^{\textup{n},\lambda})-(P_{e}^{\textup{t}})}{\R^{N_e}},
	\label{def:Delta2}
\end{equation} 
where, as before, the superscripts denote the considered behavior.

\noindent 3. We define $\Delta_3$ to be the difference between the maximum densities observed in the controlled natural and the target behavior, i.e.
\begin{equation}
	\Delta_3(\lambda):=\big|\max_{t,x}\rho^{\textup{n},\lambda}(t,x)-\max_{t,x}\rho^{\textup{t}}(t,x)\big|.
	\label{def:Delta3}
\end{equation}
Note the the maximum density is strictly related to pedestrian compression, which, in turn, is often responsible for injuries or death~\cite{helbing2007PRE}.

\smallskip

In all of the cases above we aim at finding a minimizer $\lambda^*$ such that 
\begin{equation}
	\Delta(\lambda^*)=\min_{\lambda\in\Lambda}\Delta(\lambda)
	\label{def:lambda*}
\end{equation}
along with the corresponding optimally controlled natural behavior $\rho^{\textup{n},\lambda^*}(t,x)$.

\begin{remark}
Let us stress that the result of any environment optimization is necessarily affected by both the domain $\Omega$ itself and the initial crowd distribution $\rho_0$. This is possibly the main limitation of the proposed crowd control approach. However, we believe that when the walking freedom is limited, as in structured environments like theaters or stadiums, controlling crowds by mean of environmental design can provide an effective improvement for the cost criteria $\Delta$ described above. 
\end{remark}

\subsection{Minimum search}
\label{sec:minimum_search}
Solving problem~\eqref{def:lambda*} is a hard task. Analytical methods based on the computation of shape derivatives are not applicable here, due to the complexity of the problem and to the limited stability results available in the context of multidimensional conservation laws. Moreover, as we will see in Section~\ref{sec:numerical_tests}, the environmental value function $\Delta$ is far from being differentiable, therefore methods based on the (numerical) evaluation of its gradient should be discarded. Among the derivative-free methods, we consider two options.

\paragraph{Exhaustive method} 
We define $\bar\Lambda$ to be a discrete finite approximation of the set of admissible parameters $\Lambda$ and we simply evaluate the function $\Delta$ for any $\lambda\in\bar\Lambda$. This approach can be pursued only if the dimension of $\Lambda$ is low, say at most $2$. In Section~\ref{sec:numerical_tests} we will consider the case of a single controlled obstacle with rectangular shape and fixed sides, described by means of the coordinates of its barycenter, which we denote by $(x_\mathcal{O},y_\mathcal{O})$. This way we have $\Lambda\subset\Omega\subset\R^2$. Then we compute $\Delta$ on the same grid used for approximating the density $\rho$, except for some additional constraints related to the \emph{admissibility} of the obstacle. Namely, we impose that the controlled obstacle can neither overlap the initial pedestrian distribution $\rho_0$ nor overhang the boundaries of the domain. We also require it to not completely clog the entrance or the exit, if only one of them is present.

\paragraph{Modified compass search}
Starting from some initial guess $\lambda_0$, this method generates a sequence $\{\lambda_1,\,\lambda_2,\,\ldots\}$ converging to a (local) minimum of $\Delta$. The controlled obstacle $\mathcal{O}_\lambda$ is moved in consequence of small random variations of the current value of $\lambda$. If the variation is advantageous, i.e., if $\Delta$ decreases, the variation is kept, otherwise it is discarded. In Section~\ref{sec:numerical_tests} we will consider again the case of a single rectangular obstacle, but with variable sides. Namely, the control parameters are given by the barycenter of the rectangle plus the length of its two sides, hence $\Lambda\subset\Omega\times\R_+^2\subseteq\R^4$. As a consequence, variations of the parameter $\lambda$ imply displacements of the obstacle in $\Omega$ and/or changes of its size. The same admissibility constraints as above are imposed on the obstacles ${\cal O}_\lambda$. In order to select a new value of $\lambda$, we first extract a random integer $p\in\{1,\,2,\,3,\,4,\,5\}$, then we randomly apply one of the following $8$ rules:
\begin{enumerate}
\item[{R1.--4.}] Move the obstacle $p$ cells rightward/leftward/upward/downward;
\item[{R5.--6.}] Stretch/shrink the obstacle of $2p$ cells horizontally;
\item[{R7.--8.}] Stretch/shrink the obstacle of $2p$ cells vertically;
\end{enumerate}
On the whole, we consider $5\times 8=40$ rules. Furthermore, in order to avoid remaining trapped in local minima, we complement the above rules with the simulated annealing technique. 
%
%

\section{Numerical tests}
\label{sec:numerical_tests}
In this section we compare the evacuation (minimum time) problem for basic, rational, and highly rational behaviors by means of some numerical experiments. Then, we face the challenging shape optimization problem.

Note that all parameters used in the simulations are meant to be realistic but exploratory. In particular, they do not stem from any real measurement. Moreover, we consider here crowded but not panic situations.

In all figures, exits are colored in red, entrances (if any) in black, fixed obstacles (i.e., obstacles possibly already present in the domain independently of the control/optimization procedure) in dark gray, and controlled obstacles (i.e., obstacles purposely added for control) in yellow. Pedestrian density ranges from light green, corresponding to low crowding, to dark blue, corresponding to high crowding.

\subsection{Basic, rational, and highly rational behaviors}
In order to compare the basic and rational behaviors, we consider a square room with ten exits on the upper side separated by pillars. The area of the domain is $50\times 50\unit{m^2}$. We set $\alpha=170^\circ$, $R=1.5\unit{m}$, and $F=8\unit{m^2/s}$. At the initial time, a density corresponding to $N_P=43$ pedestrians is in the room.

\begin{figure}[!t]
\begin{center}
\includegraphics[width=\textwidth]{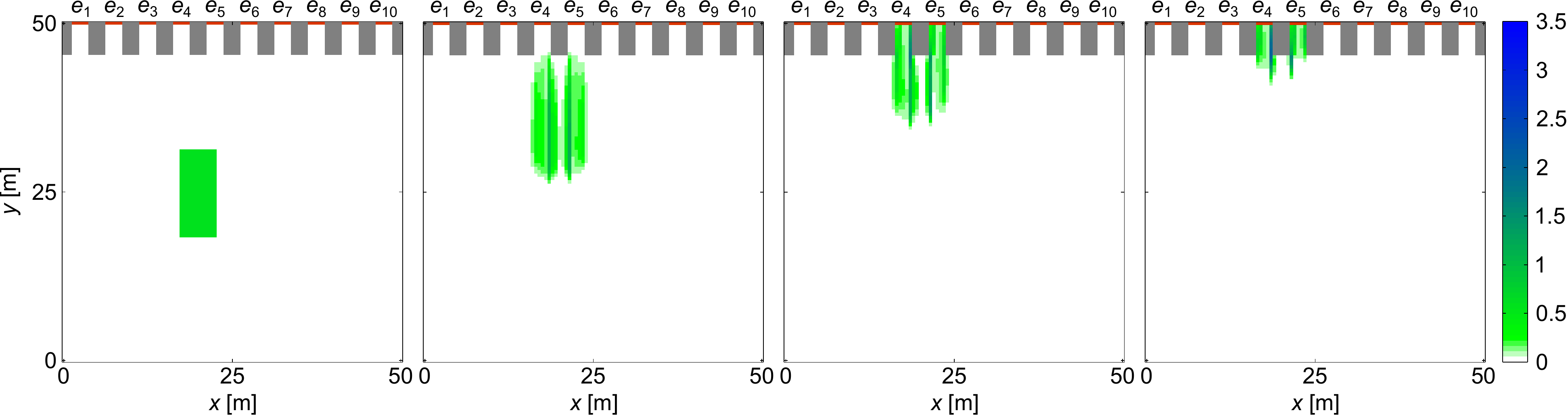}
\end{center}
\caption{Density evolution in case of basic behavior.}
\label{fig:colon_irraz_screenshots}
\end{figure}

Figure~\ref{fig:colon_irraz_screenshots} shows the basic behavior of the crowd. The group, initially arranged in a rectangular formation, points toward the nearest two exits (namely, $e_4$ and $e_5$), regardless of the fact that this choice will surely cause a congestion near them. 

\begin{figure}[!t]
\begin{center}
\includegraphics[width=\textwidth]{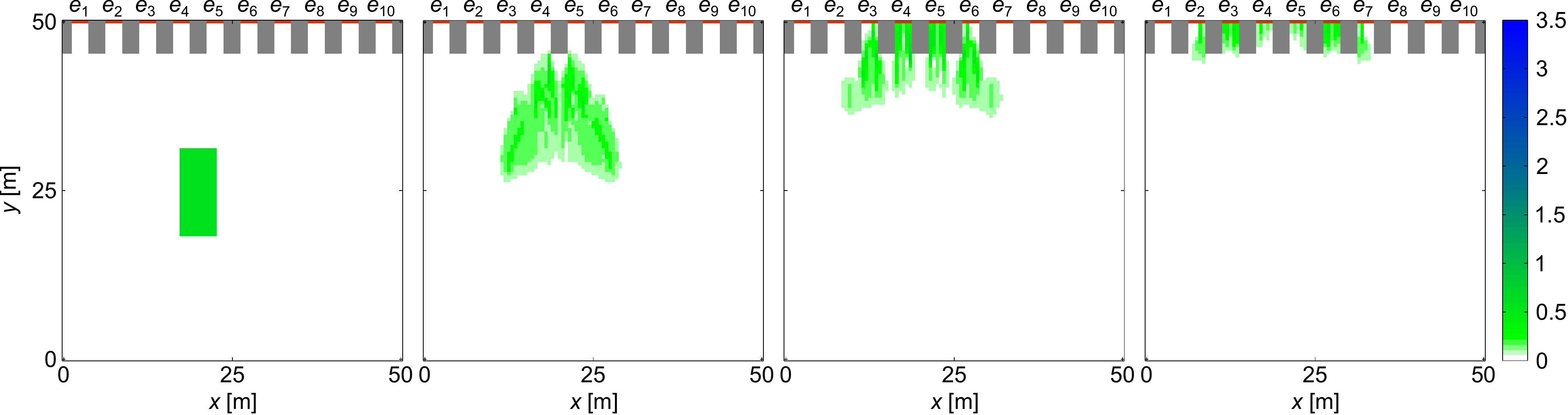}
\end{center}
\caption{Density evolution in case of rational behavior.}
\label{fig:colon_raz_screenshots}
\end{figure}

Figure~\ref{fig:colon_raz_screenshots} shows instead the rational behavior of the crowd. In this case the group spreads out immediately, because people tend to avoid congested paths. After a while, two groups of people find it convenient to reach the exits $e_3$ and $e_6$, the way to them being completely clear. Once exits $e_3$ and $e_6$ become congested, two other groups split, that point toward exits $e_2$ and $e_7$.

Three main differences are worth being noted between the two behaviors, see also Table~\ref{tab:colon_cfr}:
First, basically behaved pedestrians use two exits, rationally behaved ones use six exits.
Second, the evacuation time slightly decreases in the rational case, because the splitting of the crowd turns out to be advantageous.
Third, the maximal density is higher in the basic case, because of the congestion at exits $e_4$ and $e_5$.
\begin{table}[!t]
\caption{Comparison between basic and rational behavior}
\label{tab:colon_cfr}
\begin{center}
\begin{tabular}{|l|c|c|c|}
\hline 
Behavior & Used exits & $\tevac\unit{[s]}$ & $\rhomax\unit{[ped/m^2]}$\\
\hline\hline
Basic & 2 & 44.55 & 1.93 \\
\hline
Rational & 6 & 40.95 & 0.80 \\
\hline
\end{tabular}
\end{center}
\end{table}

\medskip

In order to compare the rational and highly rational behaviors we consider a room with no obstacles and two exits on the lower wall. The area of the domain is $50\times 25\unit{m^2}$. We set $\alpha=170^\circ$, $R=1.5\unit{m}$, and $F=8\unit{m^2/s}$. At the initial time two groups of pedestrians are in the room, precisely Group 1 ($N_P=15$ people) located at the bottom and Group 2 ($N_P=25$ people) located along the right wall.

\begin{figure}[!t]
\begin{center}
\includegraphics[width=\textwidth]{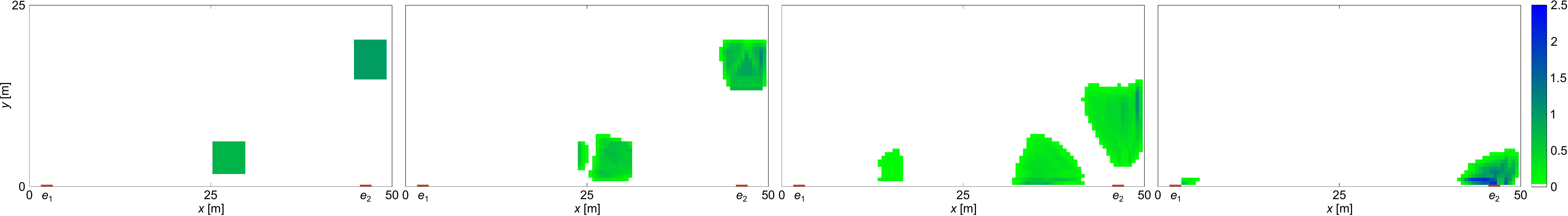}
\end{center}
\caption{Density evolution in case of rational behavior.}
\label{fig:room_usciteinbasso_raz_screenshots}
\end{figure}

Figure~\ref{fig:room_usciteinbasso_raz_screenshots} shows the rational behavior of the crowd. Most pedestrians point immediately toward the right exit, which is the nearest one. A small amount of people from Group 1 walk instead leftward, finding the way between them and the right exit too congested, hence not convenient. In this case people are not assumed to be able to foresee the evolution of the crowd at future times, therefore they cannot predict the congestion which will inevitably form near the right exit.

\begin{figure}[!t]
\begin{center}
\includegraphics[width=\textwidth]{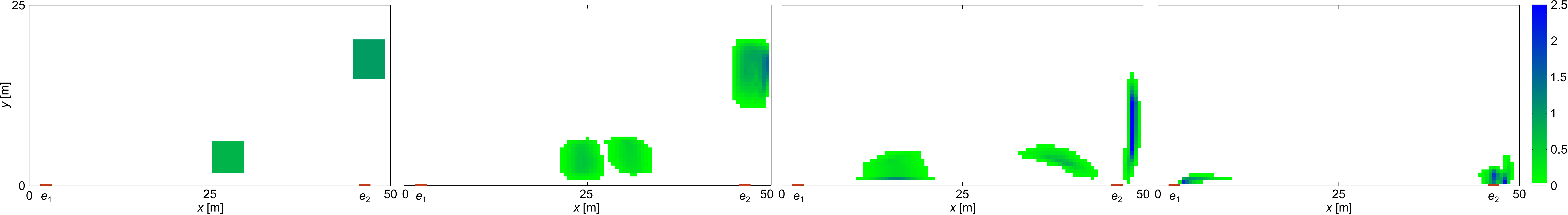}
\end{center}
\caption{Density evolution in case of highly rational behavior.}
\label{fig:room_usciteinbasso_hraz_screenshots}
\end{figure}

Figure~\ref{fig:room_usciteinbasso_hraz_screenshots} shows instead the highly rational behavior of the crowd. Group 1 splits immediately, about half of it pointing leftward. In this case people foresee the future congestion at the right exit and take countermeasures. Note also the different behaviors of the two subgroups of Group 1. The one moving to the left naturally points toward the left exit, while the other one reaches the right exit from above, sharing the exit with Group 2. The latter, in turn, reaches the exit flattened to the wall. In conclusions, in this case the flows through the two exits are quite balanced. 

An improvement in the evacuation time is also observed: $32.5\unit{s}$ for highly rational pedestrians and $36\unit{s}$ for rational pedestrians.

\subsection{Environment optimization}
Here we consider again the evacuation problem and, from now on, we assume that the natural behavior coincides with the basic one, while the target behavior coincides with the rational one. Of course, other combinations are possible. Nevertheless, we believe that considering these two behaviors is well representative of the main modeling ideas of the paper while allowing us to avoid excessively demanding numerical computations.

\subsubsection{Room with a fixed obstacle}
In this test we consider a square room with one entrance, two exits, and one fixed obstacle. The overall area of the room is $50\times 50\unit{m^2}$, while the dimension of the obstacle is $7.5\times 17\unit{m^2}$. We set $\alpha=170^\circ$, $R=1.5\unit{m}$, and $F=8\unit{m^2/s}$. Pedestrians enter the room continuously for $25\unit{s}$ at a rate of $3.5\unit{ped/s}$.

We first compute the natural and target behaviors, then we test the effect on the dynamics of an additional controlled obstacle with fixed shape, running an exhaustive optimum search. Finally, we use the compass search to find the best shape of the controlled obstacle.

\begin{figure}[!t]
\begin{center}
\includegraphics[width=.8\textwidth]{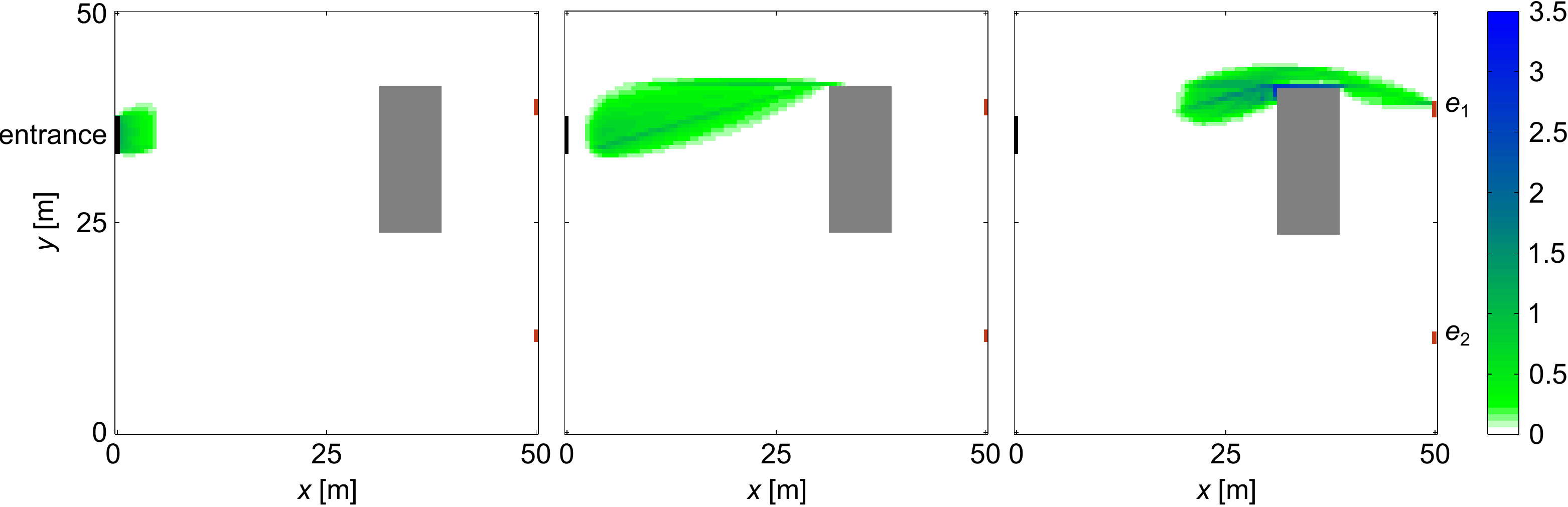}
\end{center}
\caption{Density evolution in case of natural behavior.}
\label{fig:base_ost_irraz_screenshots}
\end{figure}

Figure~\ref{fig:base_ost_irraz_screenshots} shows the natural behavior of the crowd. People point toward the upper exit $e_1$ on the right wall of the domain. While passing by the fixed obstacle pedestrians slow down a bit, because their mutual repulsion pushes some of them against the wall.

\begin{figure}[!t]
\begin{center}
\includegraphics[width=.8\textwidth]{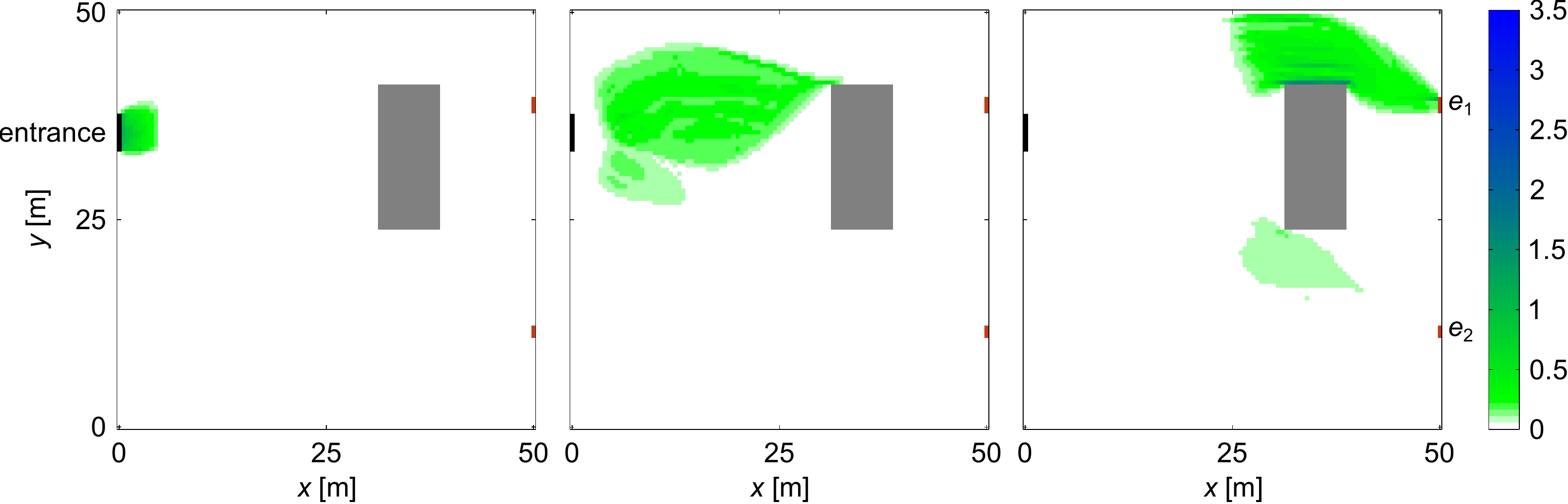}
\end{center}
\caption{Density evolution in case of target behavior.}
\label{fig:base_ost_raz_screenshots}
\end{figure}

\begin{figure}[!t]
\begin{center}
\includegraphics[width=.75\textwidth]{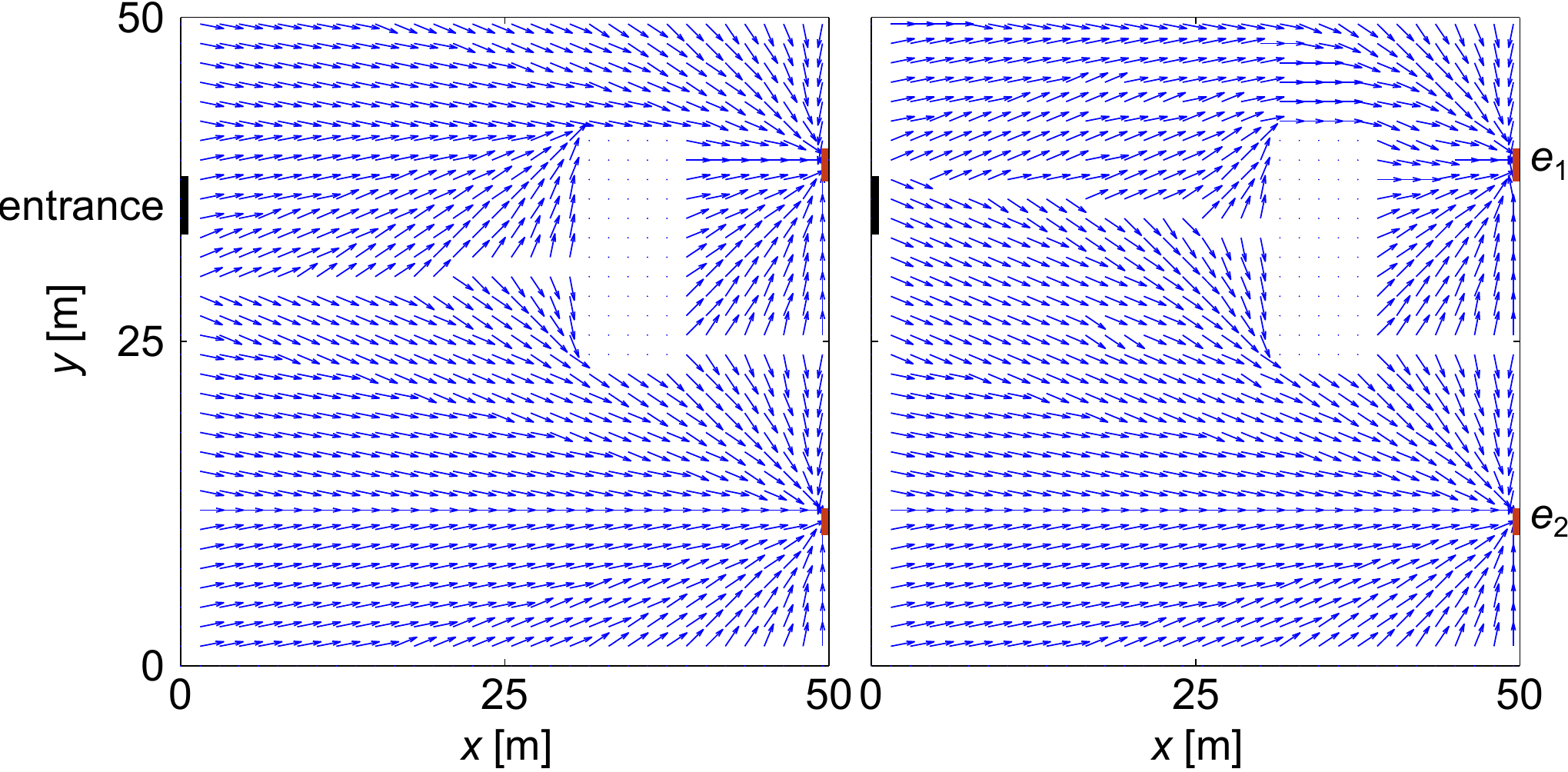}
\end{center}
\caption{Velocity $v_b^\ast$ for natural behavior (left) and for target behavior at a later time (right).}
\label{fig:base_ost_raz_Vdes}
\end{figure}

Figure~\ref{fig:base_ost_raz_screenshots} shows instead the target behavior of the crowd. Similarly to the case of the natural behavior, pedestrians point initially to the upper exit $e_1$. After a while, however, the group splits, because some individuals find it more convenient to reach the lower exit $e_2$, the way to it being completely clear. In conclusion, the space occupancy turns out to be more balanced.
The group splitting is even better understood looking at the behavioral velocity field $v_b^\ast$. Figure~\ref{fig:base_ost_raz_Vdes} shows both $v_b^\ast$ for the natural behavior (which coincides with $v_b^\ast$ for the target behavior at the initial time) and $v_b^\ast$ for the target behavior at a later time. A  discontinuity line in the velocity field is clearly visible, which coincides with the curve of nondifferentiability of the minimum time function $\phi$. Such a line divides the domain in two regions: Above the line it is more convenient to point to the upper exit $e_1$, below the line it is instead more convenient to point to the lower exit $e_2$. Along the discontinuity line either optimal path toward $e_1$ or $e_2$ can be indifferently chosen.

The differences between the natural and target behaviors are summarized in the first two rows of Table~\ref{tab:base_ost_cfr}. In this test we observe a larger difference between the two evacuation times, due to the fact that the congestion forming near the upper wall of the obstacle slows down considerably the naturally behaved pedestrians. The maximal density is also different, being lower in the target case for the same reason.
\begin{table}[!t]
\caption{Comparison between pedestrian behaviors.}
\label{tab:base_ost_cfr}
\begin{center}
\begin{tabular}{|l|c|c|c|}
\hline 
Behavior & Used exits & $\tevac\unit{[s]}$ & $\rhomax\unit{[ped/m^2]}$ \\
\hline\hline
Natural & 1 & 120.60 & 3.35 \\
\hline
Target & 2 & 95.85 & 2.37 \\
\hline
Controlled natural & 2 & 99.00 & 2.24 \\
\hline
\end{tabular}
\end{center}
\end{table}

Let us now try to make pedestrians behave more rationally than they usually do. We choose the evacuation time as the cost criterion for the environment optimization, i.e., the environmental cost function is $\Delta_1$ defined in~\eqref{def:Delta1}. We consider first the exhaustive method described in Section~\ref{sec:minimum_search}, with a square-shaped controlled obstacle of fixed area $5\times 5\unit{m^2}$, and we denote the coordinates of its barycenter by $(x_\mathcal{O},y_\mathcal{O})$. Then, for any \emph{admissible} position of this 
obstacle we compute the controlled natural evacuation time, see~\eqref{def:tevac_natural}.

\begin{figure}[!t]
\begin{center}
\includegraphics[width=\textwidth]{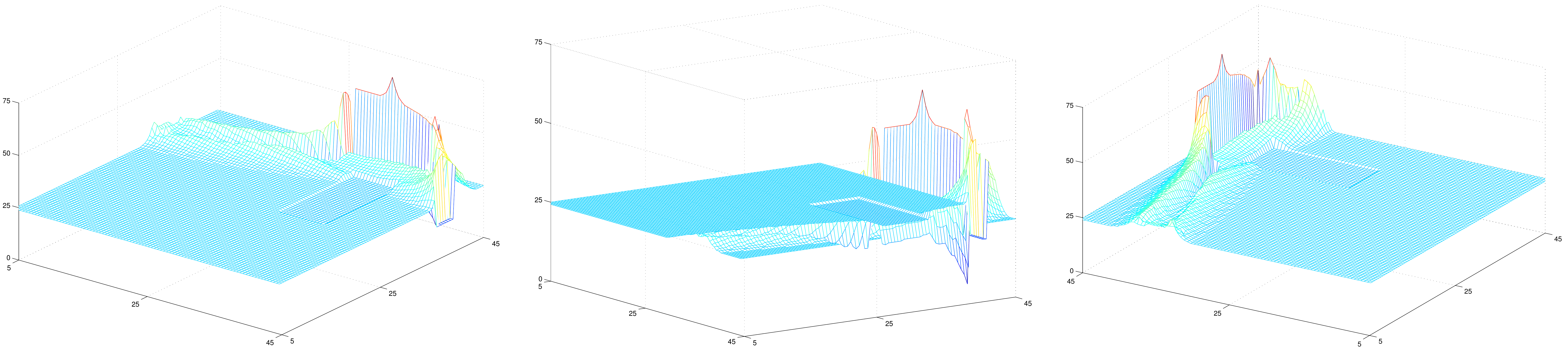}
\end{center}
\caption{The function $\Delta_1(x_\mathcal{O},y_\mathcal{O})$. From above, exit view (left), from below, exit view (center), and from above, entrance view (right).}
\label{fig:base_ost_FV}
\end{figure}

Figure~\ref{fig:base_ost_FV} shows the function $\Delta_1(x_\mathcal{O},y_\mathcal{O})$, which gives the distance between the controlled natural behavior and the target one. We immediately see that $\Delta_1$ is constant (and equal to $t^\textup{n}_\textup{evac}-t^\textup{t}_\textup{evac}$=120.60 s - 95.85 s = 24.75 s, cf.\ Table \ref{tab:base_ost_cfr}) in a large part of the domain, since there the obstacle is noninfluential. For instance, obstacle locations close to the bottom wall of the room are well outside the crowd path. Conversely, locations along the path joining the entrance and the upper exit $e_1$ heavily affect the dynamics, see the ``grand canyon''. On the top of the fixed obstacle, a high ``wall'' is visible, meaning that pedestrians are notably slowed down by obstacles in those positions, because only narrow passages remain open for walking. By looking at $\Delta_1$ from below, we clearly see one local minimum (to the left of the fixed obstacle) and one global minimum (to the right of the fixed obstacle).

\begin{figure}[!t]
\begin{center}
\includegraphics[width=0.65\textwidth]{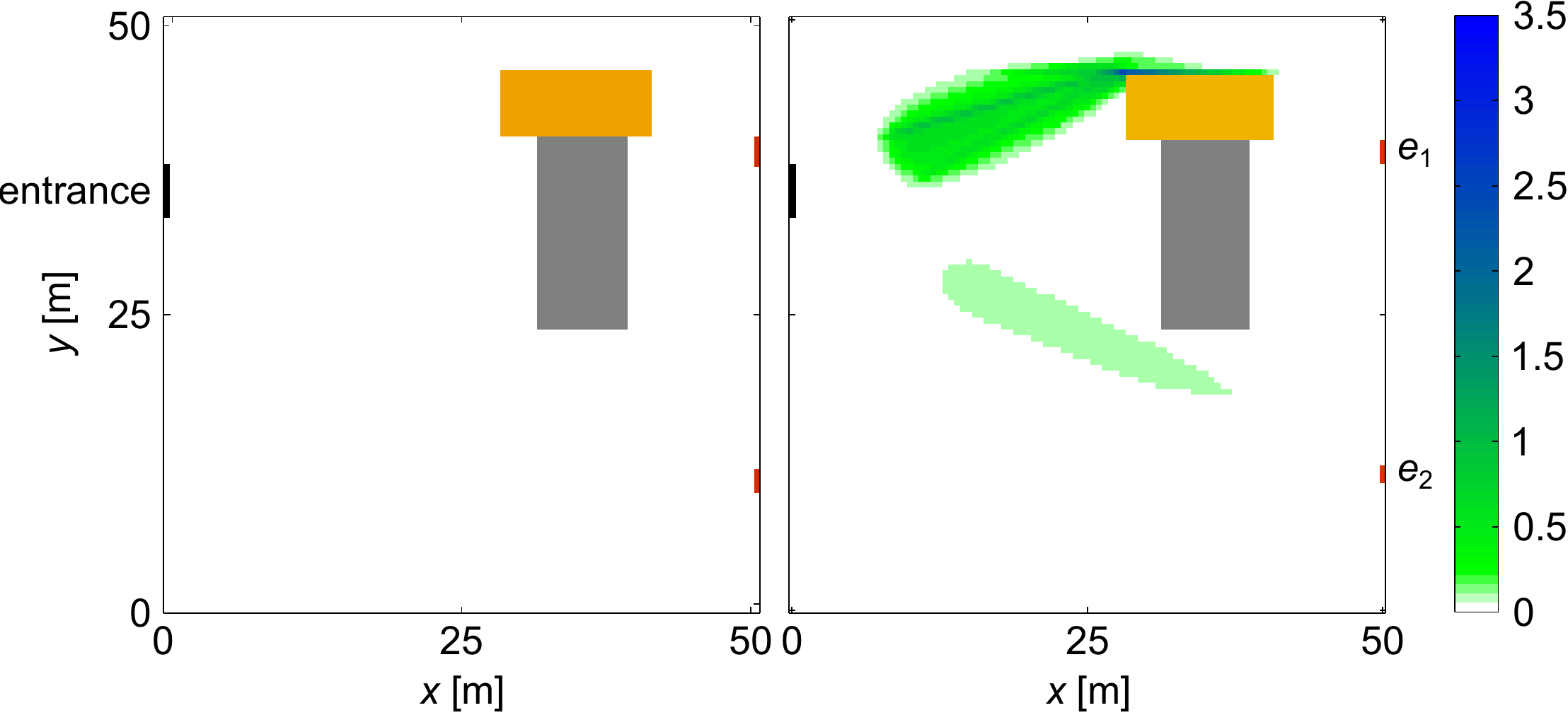}
\end{center}
\caption{A manually-tuned initial guess for the controlled obstacle (left), and the corresponding natural behavior (right).}
\label{fig:base_ost_gi_manuale}
\end{figure}

We now run the compass search method described in Section~\ref{sec:minimum_search}. Inspired by the function $\Delta_1$, we choose the controlled obstacle depicted in Fig.~\ref{fig:base_ost_gi_manuale}-left as initial guess for the descent. It is close to the global minimum but not exactly on it, and it is rectangular. The evacuation time in case of natural behavior with this additional obstacle is $125.55\unit{s}$, higher than the purely natural one, see Table~\ref{tab:base_ost_cfr}. As shown in Fig.~\ref{fig:base_ost_gi_manuale}-right, such a controlled obstacle is able to split the group of pedestrians but does not improve the evacuation time, because of the high congestion which forms on the top of the obstacle and at the upper exit. 

\begin{figure}[!t]
\begin{center}
\includegraphics[width=0.65\textwidth]{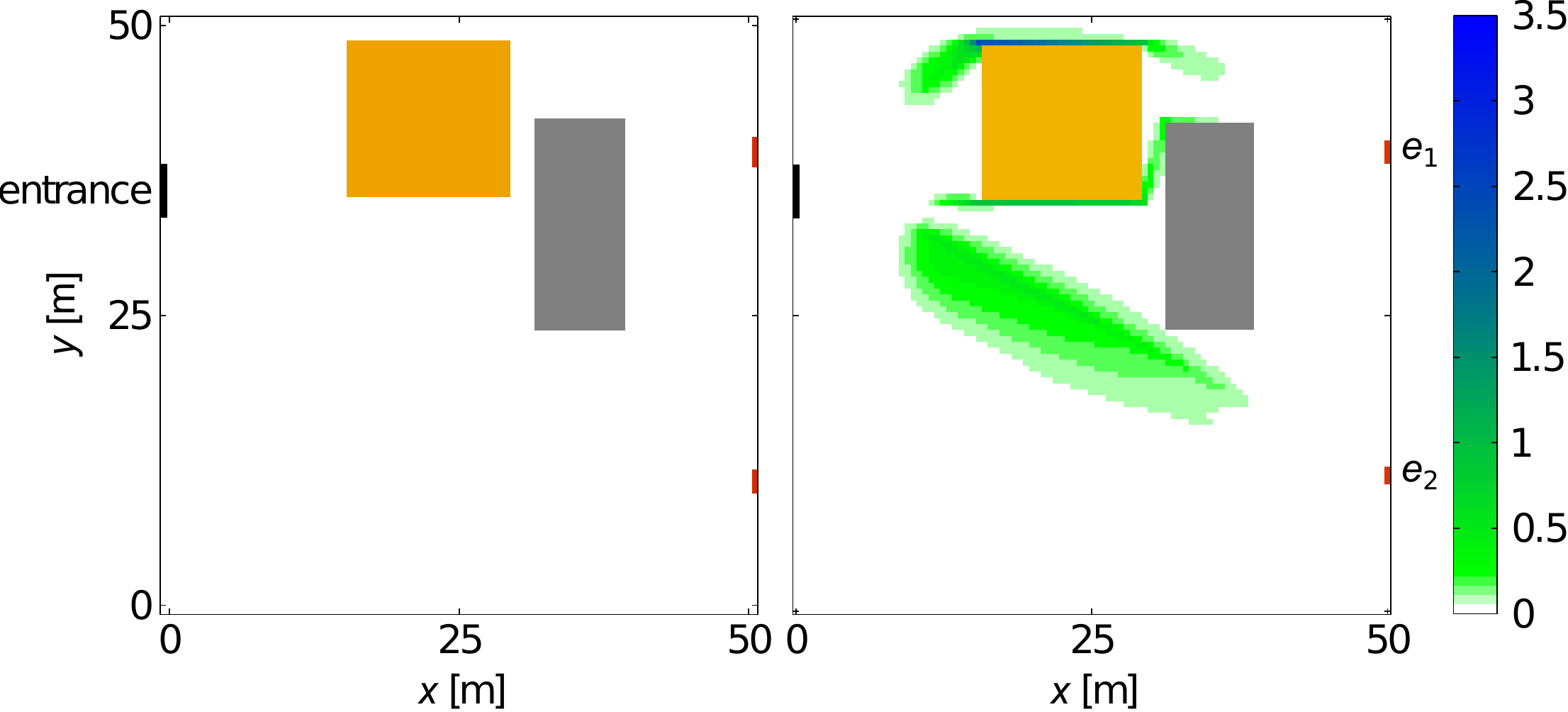}
\end{center}
\caption{Optimized controlled obstacle resulting from the compass search method (left), and the corresponding natural behavior (right).}
\label{fig:base_ost_discesa_manuale}
\end{figure}

Starting from such an initial guess, the compass search converges to a new controlled obstacle which is quite different in both shape and size, see Fig.~\ref{fig:base_ost_discesa_manuale}-left. The position of the optimized obstacle is strategical: the obstacle turns out to be able to split the group in three parts, in such a way that pedestrians use both exits $e_1$, $e_2$ giving rise to a lower congestion in front of them. A screenshot of the dynamics is shown in Fig.~\ref{fig:base_ost_discesa_manuale}-right. The evacuation time in case of such a controlled natural behavior is $99\unit{s}$ (cf. the third row of Table~\ref{tab:base_ost_cfr}), which is \emph{rather close} to the target one. This is one of the most interesting results of the paper which also reproduces, as a byproduct, the Braess' paradox.

%

\subsubsection{``Sapienza'' campus}
In this test we consider a simplified reproduction of the campus area of ``Sapienza'' University of Rome with seven exits. The area of the domain is $450\times 450\unit{m^2}$. We set $\alpha=170^\circ$, $R=13.5\unit{m}$, and $F=12\unit{m^2/s}$. At the initial time, $N_P=3500$ pedestrians are in the campus.

\begin{figure}[!t]
\begin{center}
\href{http://www.emiliano.cristiani.name/attach/Sapienza_paper_irraz.avi}
{\includegraphics[width=.8\textwidth]{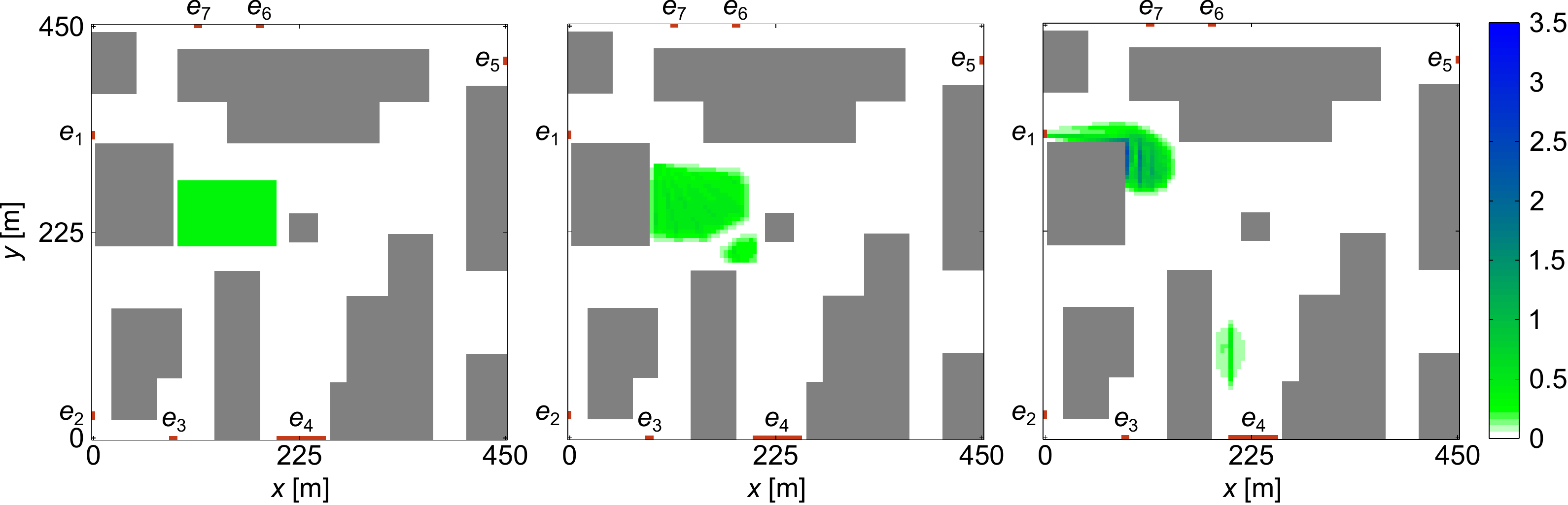}}
\end{center}
\caption{Density evolution in case of natural behavior. Click on the image to see the simulation.}
\label{fig:sap_irraz_screenshots}
\end{figure}

Figure~\ref{fig:sap_irraz_screenshots} shows the natural behavior of the crowd in the predefined structural arrangement of the campus. The group immediately splits in two parts, forming an important congestion around one corner of the building placed near exit $e_1$. At the end, people actually use only exits $e_1$ and $e_4$ for their evacuation.

\begin{figure}[t!]
\begin{center}
\href{http://www.emiliano.cristiani.name/attach/Sapienza_paper_raz.avi}
{\includegraphics[width=.8\textwidth]{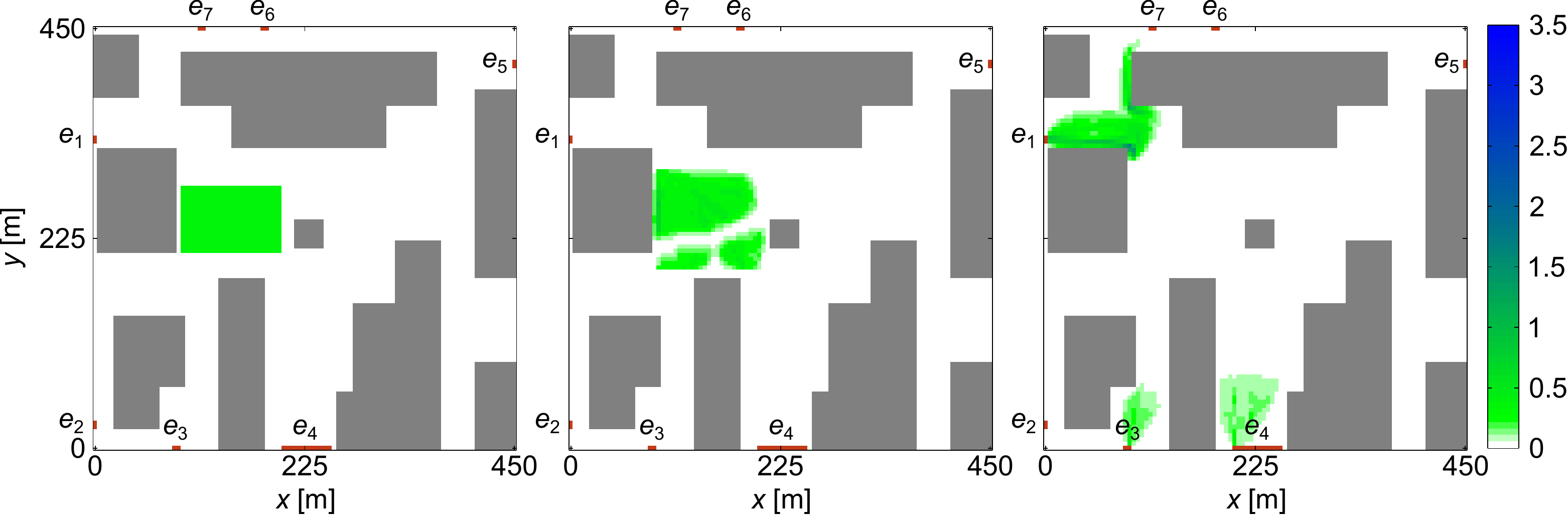}}
\end{center}
\caption{Density evolution in case of target behavior. Click on the image to see the simulation.}
\label{fig:sap_raz_screenshots}
\end{figure}

Figure~\ref{fig:sap_raz_screenshots} shows instead the target behavior of the crowd. This time pedestrians split in three groups, since people in the bottom part of the initial formation are able to estimate the slowing effect produced by the rest of the mass. After a while, the exit $e_1$ becomes congested and some individuals split again, pointing toward exit $e_7$. At the end, people actually use exits $e_1$, $e_3$, $e_4$, and $e_7$, the other exits being completely unused.
The differences between these behaviors are summarized in the first two rows of Table~\ref{tab:sap_cfr_nat_vs_target}, which also reports the percentages of pedestrians choosing each exit.
\begin{table}[t!]
\caption{Comparison between pedestrian behaviors.}
\label{tab:sap_cfr_nat_vs_target}
\begin{center}
\begin{tabular}{|l|c|c|c|l|l|l|l|}
\hline 
\multirow{2}{*}{Behavior} & Used & $\tevac$ & $\rhomax$ & \multirow{2}{*}{\%~$e_1$} & \multirow{2}{*}{\%~$e_3$} & \multirow{2}{*}{\%~$e_4$} & \multirow{2}{*}{\%~$e_7$} \\
& exits & $\unit{[s]}$ & $\unit{[ped/m^2]}$ & & & & \\
\hline\hline
Natural & 2 & 599.40 & 2.87 & 92.40 & 0.00 & 7.60 & 0.00 \\
\hline
Target & 4 & 429.30 & 1.58 & 60.95 & 9.47 & 14.88 & 14.70 \\
\hline
Controlled natural & 4 & 571.05 & 2.21 & 59.95 & 12.48 & 18.23 & 9.34 \\
\hline
\end{tabular}
\end{center}
\end{table}

\begin{figure}[!t]
\begin{center}
\includegraphics[width=0.9\textwidth]{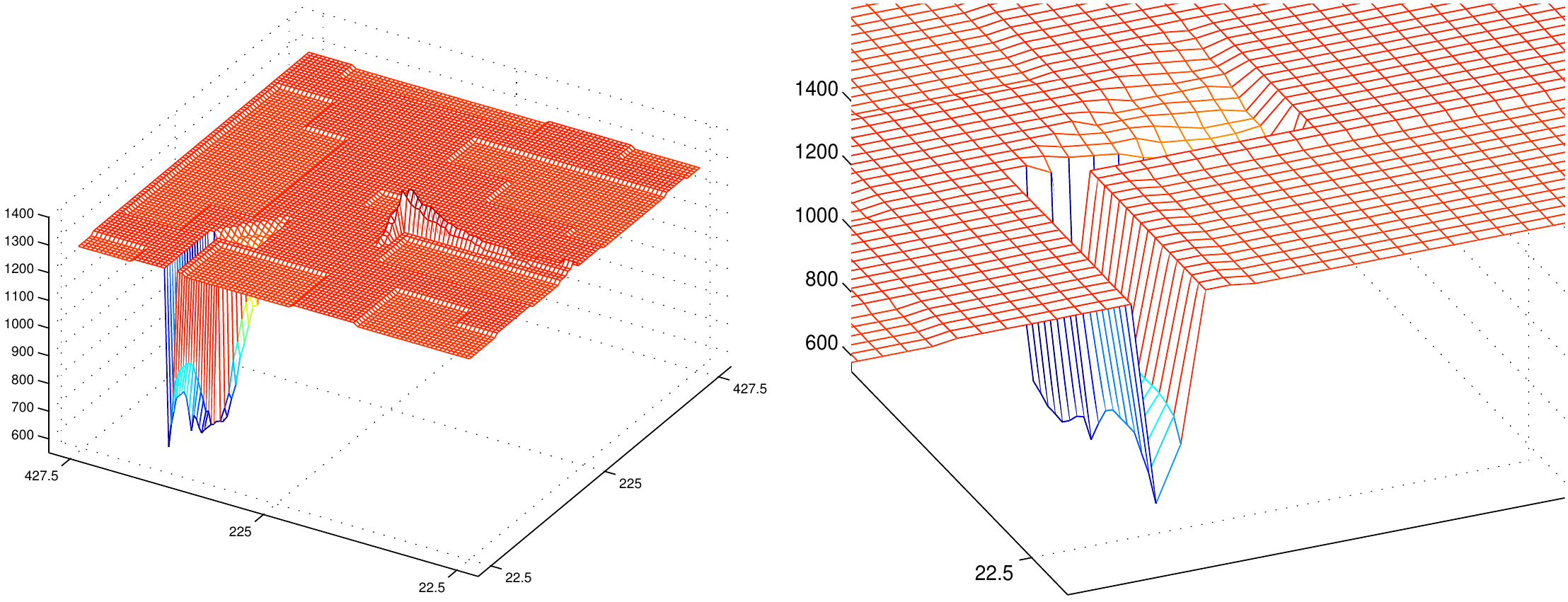}
\end{center}
\caption{The function $\Delta_2(x_\mathcal{O},y_\mathcal{O})$ (left), zoom around the downward spikes (right).}
\label{fig:sap_FV}
\end{figure}

In order to tackle the environment optimization, in this case we choose the pedestrian distribution at the exits as a cost criterion, namely the environmental cost function is $\Delta_2$ defined in~\eqref{def:Delta2}. We start with the exhaustive method, considering a square-shaped controlled obstacle of fixed area $22.5\times 22.5\unit{m^2}$. As before, we denote the coordinates of its barycenter by $(x_\mathcal{O},y_\mathcal{O})$ and for any of its admissible positions we simulate the controlled natural behavior of the crowd and compute the corresponding controlled natural cost. Figure~\ref{fig:sap_FV} shows the result of the exhaustive minimum search, i.e., the function $\Delta_2(x_\mathcal{O},y_\mathcal{O})$, which gives the distance between the controlled natural behavior and the target one. We see that $\Delta_2$ is constant in a large part of the domain, meaning that there the obstacle is noninfluential. Conversely, along the paths joining the initial crowd density and exits $e_1$, $e_4$, respectively, the obstacle affects the dynamics. In particular, a series of downward spikes are visible near exit $e_1$.

\begin{figure}[!t]
\begin{center}
\includegraphics[width=0.65\textwidth]{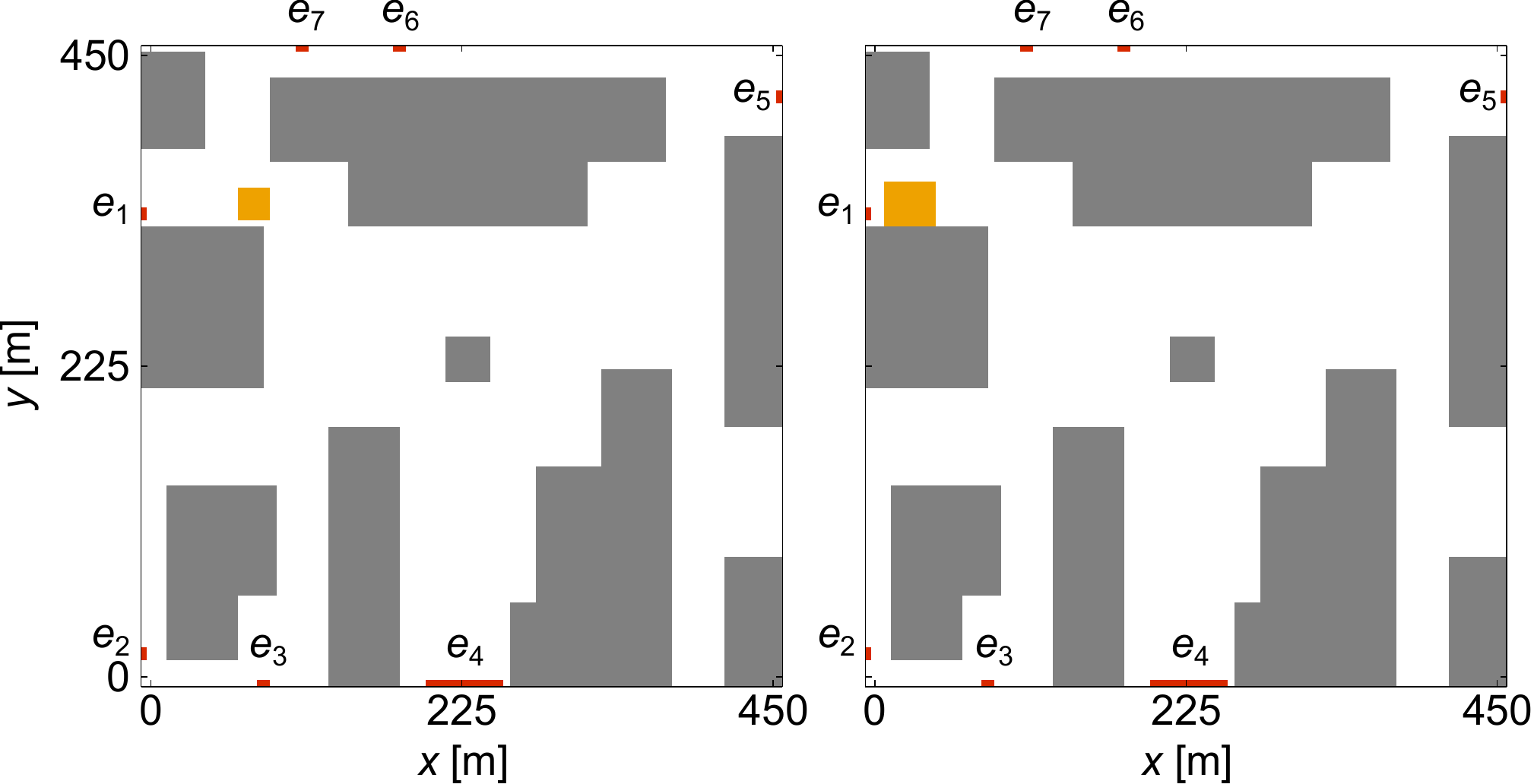}
\end{center}
\caption{Best controlled obstacle with fixed shape from the exhaustive minimum search (left), and the result of the compass search (right).}
\label{fig:sap_cfr_ostacoli_gi_esaustivo}
\end{figure}

\begin{figure}[t!]
\begin{center}
\href{http://www.emiliano.cristiani.name/attach/Sapienza_paper_irraz_ctrld.avi}
{\includegraphics[width=\textwidth]{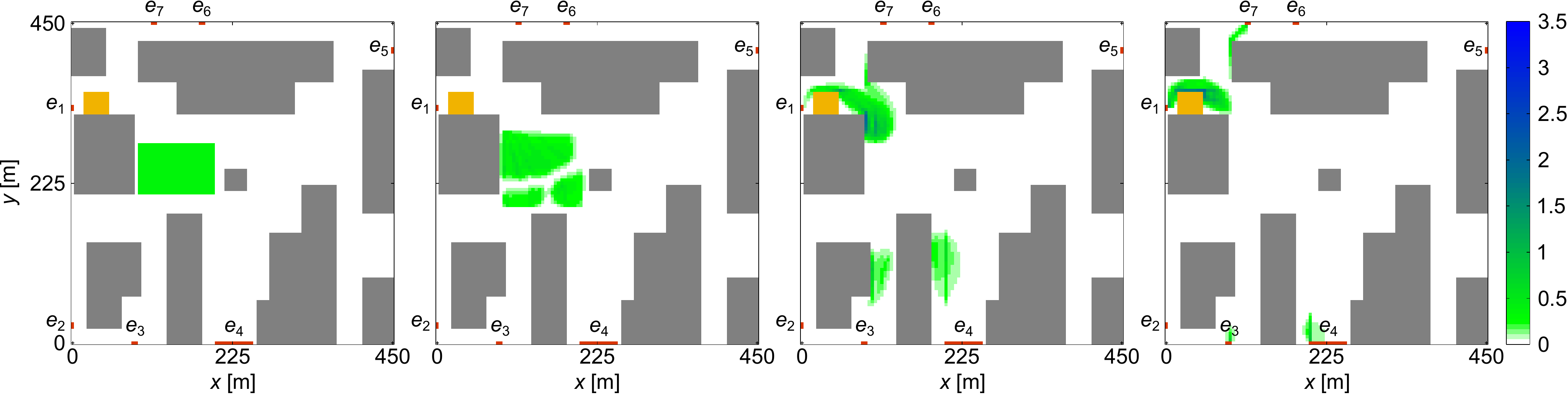}}
\end{center}
\caption{Density evolution in case of natural behavior with the best controlled obstacle resulting from the compass search. Click on the image to see the simulation.}
\label{fig:sap_gi_esaustivo_screenshots}
\end{figure}

The controlled obstacle corresponding to the global minimum of $\Delta_2$, see Fig.~\ref{fig:sap_cfr_ostacoli_gi_esaustivo}-left, serves as an initial guess for the compass search. The method converges to the optimized obstacle depicted in Fig.~\ref{fig:sap_cfr_ostacoli_gi_esaustivo}-right, which is larger than the initial guess and closer to exit $e_1$.
Such an optimal obstacle turns out to have very nice effects on the density evolution in this scenario, see Fig.~\ref{fig:sap_gi_esaustivo_screenshots}. Surprisingly, the controlled obstacle is able to split the group in \emph{four} parts (cf. the third row of Table~\ref{tab:sap_cfr_nat_vs_target}), thereby determining an effective usage of the available exits which is rather close to the one of the target behavior. Moreover, the evacuation time and the maximal density are comprised between the natural and the target case.

\section{Open problems and research perspectives}
\label{sec:openpb}
%
Each level of ``rationality'' presented in this paper leads to a different analytical framework where the corresponding model can be naturally set, but all of them present some common challenges. 

The first one comes from control theory. Indeed, models presented in Section~\ref{sec:models} require pedestrians to optimize the perceived dynamics~\eqref{eq:dyn_perc} w.r.t.\ a given performance criterion $\mathcal{C}$ (typically of the form~\eqref{eq:min_time} or~\eqref{eq:fixed_horiz}) and to construct a $\mathcal{C}$-optimal feedback $v_b^\ast$ on the whole domain $\Omega$. This is actually a nontrivial problem, because optimal feedbacks are generally neither Lipschitz continuous nor even continuous. Consider, for instance, the case of basic behavior with cost~\eqref{eq:min_time}: In this case $v_b^\ast$ is given by~\eqref{eq:basic_optimal1} but $\phi$ is usually at most Lipschitz continuous, thus $v_b^\ast$ is expected to be only in $L^\infty(\Omega)$. In turn, such a lack of regularity makes unclear the sense in which the resulting closed loop system obtained by~\eqref{eq:dyn_perc} with the feedback $v_b^\ast$, i.e.,
\begin{equation}
	\dot{y}(t)=f\big(t,y(t),v_b^\ast(y(t))\big)
	\label{eq:closed_loop}
\end{equation}
admits a solution, given that Carath\'eodory solutions do not generally exist for arbitrarily discontinuous vector fields. 

A possible way to overcome this issue is to relax the optimality requirement with \emph{nearly optimality} for the given cost $\mathcal{C}$. Indeed, in several cases one can construct $\varepsilon$-optimal feedbacks $v_b^\ast$ whose discontinuities are tame enough to guarantee the existence of Carath\'eodory solutions to~\eqref{eq:closed_loop} for every initial datum and to ensure good stability properties, see~\cite{AB99,BressanPriuli}. However, conditions have to be imposed on $\partial\Omega$ to both the vector field $f$ and the cost functions $\ell$, $g$ in order to ensure that~\eqref{eq:bc_v} holds.

A second problem is a consequence of the aforesaid lack of continuity at the level of the conservation law~\eqref{eq:conslaw}. 
When the vector field $v$ in~\eqref{eq:v} is not Lipschitz continuous, the solution to~\eqref{eq:conslaw} evolves as a measure concentrated on the whole set of trajectories of the control system~\eqref{eq:dyn_real}, which are in general not unique. Then one has to combine the techniques from~\cite{TosinFrasca} for nonlocal velocities with those from~\cite{AGSbook} for vector fields with low regularity, and take into account~\eqref{eq:bc_v}.

For highly rational behaviors, there is the additional difficulty of the coupling between the PDEs in~\eqref{eq:mfg} and~\eqref{eq:HJB_hrational_MT}. A generalization of the results in~\cite{LL} to the case treated in our model, with the presence of the additional coupling term $v_i[\rho]$, is given in~\cite{PriuliMFG}: In particular, it is shown that a solution pair $(\rho,\phi)$ exists provided the interaction velocity field $v_i$ and the cost functions $\ell$, $g$ are smooth enough. Even in this context, however, the case of a bounded domain $\Omega$ with impermeability boundary condition~\eqref{eq:bc_v} is still completely open. The uniqueness of Nash equilibria, when they exist, is open as well: standard techniques fail due to the presence of the nonlocal vector field in the HJB equation, and analogies with other differential game contexts where uniqueness fails (see~\cite[Section~1.3]{cristiani2014pp} and references therein), strongly suggest that only additional conditions on the dynamics and on the interaction terms might lead to a unique Nash equilibria.


Other forms of rationality may be also taken into account, linked e.g., to cooperation in subgroups of individuals (cf. the notion of \emph{social groups}~\cite{moussaid2010PLOS}). For instance, during evacuation people may tend to seek out friends and family members rather than acting selfishly. This clearly impacts on their egress strategy, namely on their choice of $v_b^\ast$. Modeling such an attitude would possibly require the use of multiscale techniques, such as those developed in~\cite{cristiani2011MMS,cristiani2014book}, in order to catch the partly singular behavior of few individuals with respect to the rest of the anonymous collective crowd.




Regarding research perspectives, one could apply the ideas presented above to other phenomena where the evolution of a macroscopic quantity is influenced by (optimal) choices of single individuals at microscopic level. In~\cite{cristiani2014pp}, the first and second authors consider the same rationality degrees introduced in section~\ref{sect:behavioral} in the context of traffic flow models on a network of roads. In such a situation, drivers can decide their own path at each crossroad, based on different amounts of information, and this reflects heavily into the resulting macroscopic flow of the car density.

Finally, we note that this paper promotes a way to indirectly control crowds: By using environmental controls, the natural behavior of single pedestrians is preserved while surrounding conditions are changed. Following the same broad idea, paper~\cite{albi2015arxiv} proposes, as an alternative, to ``hide'' in the crowd a few special individuals, who are suitably trained to steer the mass toward a specific target. These special individuals are not recognized as such by other people, thereby preserving the natural (unaware) behavior of the crowd. New models and experiments with real pedestrians are presented to validate this approach.


\bibliographystyle{siam}
\bibliography{main}

\begin{thebibliography}{10}

\bibitem{agnelli2015M3AS}
{\sc J.~P. Agnelli, F.~Colasuonno, and D.~Knopoff}, {\em A kinetic theory
  approach to the dynamic of crowd evacuation from bounded domains}, Math.
  Models Methods Appl. Sci., 25 (2015), pp.~109--129.

\bibitem{albi2015arxiv}
{\sc G.~Albi, M.~Bongini, E.~Cristiani, and D.~Kalise}, {\em Invisible control
  of self-organizing agents leaving unknown environments}.
\newblock Preprint arXiv:1504.04064.

\bibitem{AGSbook}
{\sc L.~Ambrosio, N.~Gigli, and G.~Savar{\'e}}, {\em Gradient flows in metric
  spaces and in the space of probability measures}, Lectures in Mathematics ETH
  Z{\"u}rich, Birkh{\"a}user Verlag, 2~ed., 2008.

\bibitem{AB99}
{\sc F.~Ancona and A.~Bressan}, {\em Patchy vector fields and asymptotic
  stabilization}, ESAIM Control Optim. Calc. Var., 4 (1999), pp.~445--471.

\bibitem{bardibook}
{\sc M.~Bardi and I.~Capuzzo~Dolcetta}, {\em Optimal control and viscosity
  solutions of {H}amilton-{J}acobi-{B}ellman equations}, Birkh\"auser, Boston,
  1997.

\bibitem{bellomo2011SR}
{\sc N.~Bellomo and C.~Dogb\'{e}}, {\em On the modeling of traffic and crowds:
  {A} survey of models, speculations, and perspectives}, SIAM Rev., 53 (2011),
  pp.~409--463.

\bibitem{bokanowski2011IFAC}
{\sc O.~Bokanowski and H.~Zidani}, {\em Minimal time problems with moving
  targets and obstacles}, in Proceedings of the 18th {IFAC} {W}orld {C}ongress,
  2011, pp.~2589--2593.

\bibitem{braess2005TS}
{\sc D.~Braess, A.~Nagurney, and T.~Wakolbinger}, {\em On a paradox of traffic
  planning}, Transport. Sci., 39 (2005), pp.~446--450.

\bibitem{BressanPriuli}
{\sc A.~Bressan and F.~S. Priuli}, {\em Nearly optimal patchy feedbacks},
  Discrete Contin. Dyn. Syst. Ser. A, 21 (2008), pp.~687--701.

\bibitem{BuDFMaWo}
{\sc M.~Burger, M.~Di~Francesco, P.~A. Markowich, and M.-T. Wolfram}, {\em Mean
  field games with nonlinear mobilities in pedestrian dynamics}, Discrete
  Contin. Dyn. Syst. Ser. B, 19 (2014), pp.~1311--1333.

\bibitem{cacace2014SISC}
{\sc S.~Cacace, E.~Cristiani, and M.~Falcone}, {\em Can local single-pass
  methods solve any stationary {H}amilton-{J}acobi-{B}ellmann equation?}, SIAM
  J. Sci. Comput., 36 (2014), pp.~A570--A587.

\bibitem{cardanotes}
{\sc P.~Cardaliaguet}, {\em Notes on mean field games}.

\bibitem{carlini2014SINUM}
{\sc E.~Carlini and F.~Silva}, {\em A fully discrete semi-{L}agrangian scheme
  for a first order mean field game problem}, SIAM J. Numer. Anal., 52 (2014),
  pp.~45--67.

\bibitem{mason2012SICON}
{\sc Y.~Chitour, F.~Jean, and P.~Mason}, {\em Optimal control models of
  goal-oriented human locomotion}, SIAM J. Control Optim., 50 (2010),
  pp.~147--170.

\bibitem{cristiani2011MMS}
{\sc E.~Cristiani, B.~Piccoli, and A.~Tosin}, {\em Multiscale modeling of
  granular flows with application to crowd dynamics}, Multiscale Model. Simul.,
  9 (2011), pp.~155--182.

\bibitem{cristiani2012CDC}
\leavevmode\vrule height 2pt depth -1.6pt width 23pt, {\em How can macroscopic
  models reveal self-organization in traffic flow?}, in Proceedings of the 51st
  {IEEE} {C}onference on {D}ecision and {C}ontrol, Maui, HI, USA, December
  2012, pp.~6989--6994.

\bibitem{cristiani2014book}
\leavevmode\vrule height 2pt depth -1.6pt width 23pt, {\em Multiscale Modeling
  of Pedestrian Dynamics}, {Modeling, Simulation \& Applications}, Springer,
  2014.

\bibitem{cristiani2014pp}
{\sc E.~Cristiani and F.~S. Priuli}, {\em A destination-preserving model for
  simulating {W}ardrop equilibria in traffic flow on networks}, To appear in
  Netw. Heterog. Media.

\bibitem{duives2013TRC}
{\sc D.~C. Duives, W.~Daamen, and S.~P. Hoogendoorn}, {\em State-of-the-art
  crowd motion simulation models}, Transportation Res. C, 37 (2013),
  pp.~193--209.

\bibitem{falconebook}
{\sc M.~Falcone and R.~Ferretti}, {\em Semi-Lagrangian approximation schemes
  for linear and Hamilton-Jacobi equations}, SIAM, 2014.

\bibitem{helbing2000N}
{\sc D.~Helbing, I.~Farkas, and T.~Vicsek}, {\em Simulating dynamical features
  of escape panic}, Nature, 407 (2000), pp.~487--490.

\bibitem{helbing2007PRE}
{\sc D.~Helbing, A.~Johansson, and H.~Z. Al-Abideen}, {\em Dynamics of crowd
  disasters: {A}n empirical study}, Phys. Rev. E, 75 (2007), pp.~046109/1--7.

\bibitem{helbing1995PRE}
{\sc D.~Helbing and P.~Moln\'{a}r}, {\em Social force model for pedestrian
  dynamics}, Phys. Rev. E, 51 (1995), pp.~4282--4286.

\bibitem{hoogendoorn2003OCAM}
{\sc S.~P. Hoogendoorn and P.~H.~L. Bovy}, {\em Simulation of pedestrian flows
  by optimal control and differential games}, Optim. Control Appl. Meth., 24
  (2003), pp.~153--172.

\bibitem{hoogendoorn2004TRBb}
\leavevmode\vrule height 2pt depth -1.6pt width 23pt, {\em Dynamic user-optimal
  assignment in continuous time and space}, Transportation Res. B, 38 (2004),
  pp.~571--592.

\bibitem{hughes2002TRB}
{\sc R.~L. Hughes}, {\em A continuum theory for the flow of pedestrians},
  Transportation Res. B, 36 (2002), pp.~507--535.

\bibitem{hughes2003ARFM}
\leavevmode\vrule height 2pt depth -1.6pt width 23pt, {\em The flow of human
  crowds}, Annu. Rev. Fluid Mech., 35 (2003), pp.~169--182.

\bibitem{johansson2007}
{\sc A.~Johansson and D.~Helbing}, {\em Pedestrian flow optimization with a
  genetic algorithm based on boolean grids}, in Pedestrian and Evacuation
  Dynamics 2005, N.~Waldau, P.~Gattermann, H.~Knoflacher, and M.~Schreckenberg,
  eds., Springer-Verlag Berlin Heidelberg, 2007, pp.~267--272.

\bibitem{kachroo2008book}
{\sc P.~Kachroo, S.~J. Al-nasur, S.~A. Wadoo, and A.~Shende}, {\em Pedestrian
  dynamics. Feedback control of crowd evacuation}, Understanding Complex
  Systems, Springer-Verlag, Berlin Heidelberg, 2008.

\bibitem{lachapelle2011TRB}
{\sc A.~Lachapelle and M.-T. Wolfram}, {\em On a mean field game approach
  modeling congestion and aversion in pedestrian crowds}, Transportation Res.
  B, 45 (2011), pp.~1572--1589.

\bibitem{LL}
{\sc J.-M. Lasry and P.-L. Lions}, {\em Mean field games}, Jpn. J. Math., 2
  (2007), pp.~229--260.

\bibitem{moussaid2010PLOS}
{\sc M.~Moussa\"id, N.~Perozo, S.~Garnier, D.~Helbing, and G.~Theraulaz}, {\em
  The walking behaviour of pedestrian social groups and its impact on crowd
  dynamics}, PLoS One, 5 (2010), p.~e10047.

\bibitem{piccoli2013AAM}
{\sc B.~Piccoli and F.~Rossi}, {\em Transport equation with nonlocal velocity
  in {W}asserstein spaces: convergence of numerical schemes}, Acta Appl. Math.,
  124 (2013), pp.~73--105.

\bibitem{piccoli2011ARMA}
{\sc B.~Piccoli and A.~Tosin}, {\em Time-evolving measures and macroscopic
  modeling of pedestrian flow}, Arch. Ration. Mech. Anal., 199 (2011),
  pp.~707--738.

\bibitem{PriuliMFG}
{\sc F.~S. Priuli}, {\em First order mean field games in crowd dynamics}.
\newblock Preprint arXiv:1402.7296, 2014.

\bibitem{robin2009TRB}
{\sc T.~Robin, G.~Antonini, M.~Bierlaire, and J.~Cruz}, {\em Specification,
  estimation and validation of a pedestrian walking behavior model},
  Transportation Res. B, 43 (2009), pp.~36--56.

\bibitem{shukla2009}
{\sc P.~K. Shukla}, {\em Genetically optimized architectural designs for
  control of pedestrian crowds}, in Artificial life: Borrowing from biology,
  K.~Korb, M.~Randall, and T.~Hendtlass, eds., vol.~5865 of LNCS,
  Springer-Verlag Berlin Heidelberg, 2009, pp.~22--31.

\bibitem{TosinFrasca}
{\sc A.~Tosin and P.~Frasca}, {\em Existence and approximation of probability
  measure solutions to models of collective behaviors}, Netw. Heterog. Media, 6
  (2011), pp.~561--596.

\bibitem{twarogowska2014AMM}
{\sc M.~Twarogowska, P.~Goatin, and R.~Duvigneau}, {\em Macroscopic modeling
  and simulations of room evacuation}, Appl. Math. Model., 38 (2014),
  pp.~5781--5795.

\bibitem{xia2009PRE}
{\sc Y.~Xia, S.~C. Wong, and C.-W. Shu}, {\em Dynamic continuum pedestrian flow
  model with memory effect}, Phys. Rev. E, 79 (2009), pp.~066113/1--8.

\end{thebibliography}
\end{document}